\documentclass[11pt,oneside]{article}
\usepackage{amsfonts}
\usepackage{amsmath}
\usepackage{amsthm}
\usepackage[mathscr]{euscript}
\usepackage{mathrsfs}
\usepackage{indentfirst}
\usepackage{xcolor}

\newcommand{\lbl}[1]{\label{#1}}

\newtheorem{theo}{Theorem}[section]
\newtheorem{prop}{Proposition}[section]
\newtheorem{lem}{Lemma}[section]

\newtheorem{defi}{Definition}[section]
\newcommand{\be}{\begin{equation}}
\newcommand{\ee}{\end{equation}}
\newcommand\bes{\begin{eqnarray}} \newcommand\ees{\end{eqnarray}}
\newcommand{\bess}{\begin{eqnarray*}}
\newcommand{\eess}{\end{eqnarray*}}

\newcommand\ep{\varepsilon}
\newcommand\kk{\left}
\newcommand\rr{\right}
\newcommand\dd{\displaystyle}
\newcommand\vp{\varphi}

\newcommand\ap{\alpha}

\newcommand\yy{\infty}

\newcommand\ud{\underline}
\def\theequation{\arabic{section}.\arabic{equation}}
\begin{document}
\setlength{\baselineskip}{16pt} \pagestyle{myheadings}

\begin{center}{\Large\bf Dynamics for a diffusive prey-predator model}\\[2mm]
{\Large\bf with different free boundaries}\footnote{This work was
supported by NSFC Grant 11371113}\\[4mm]
 {\Large  Mingxin Wang\footnote{{\sl E-mail}: mxwang@hit.edu.cn}}\\[0.5mm]
{\small Department of Mathematics, Harbin Institute of Technology, Harbin 150001, PR}
\\[2mm]
{\Large Yang Zhang}\\[1mm]
{\small Department of Mathematics, Harbin Engineering University, Harbin 150001, PR China}
\end{center}

\begin{quote}
\noindent{\bf Abstract.} To understand the spreading and interaction of prey and predator, in this paper we study the dynamics of the diffusive Lotka-Volterra type prey-predator model with different free boundaries. These two free boundaries, which may intersect each other as time evolves, are used to describe the spreading of prey and predator. We investigate the existence and uniqueness, regularity and uniform estimates, and long time behaviors of global solution. Some sufficient conditions for spreading and vanishing are established. When spreading occurs, we provide the more accurate limits of $(u,v)$ as $t\to\infty$, and give some estimates of asymptotic spreading speeds of $u,v$ and asymptotic speeds of $g,h$. Some realistic and significant spreading phenomena are found.


 \end{quote}

 \section{Introduction}
 \setcounter{equation}{0} {\setlength\arraycolsep{2pt}

The spreading and vanishing of multiple species is an important content in understanding ecological complexity. In order to study the spreading and vanishing phenomenon, many mathematical models have been established. In this paper we consider the diffusive Lotka-Volterra type prey-predator model with different free boundaries. It is a meaningful subject, because the following phenomenon will happen constantly in the real world:

There is some kind of species (the indigenous species, prey $u$) in a bounded area (initial habitat, for example, $\Omega_0$), and at some time (initial time, $t=0$) another kind of species (the new or invasive species, predator $v$) enters a part $\Sigma_0$ of $\Omega_0$.

In general, both species have tendencies to emigrate from boundaries to obtain their respective new habitats. That is, as time $t$ increases, $\Omega_0$ and $\Sigma_0$ will evolve into expanding regions $\Omega(t)$ and $\Sigma(t)$ with expanding fronts $\partial\Omega(t)$ and $\partial\Sigma(t)$, respectively. The initial functions $u_0(x)$ and $v_0(x)$ will evolve into positive functions $u(t,x)$ and $v(t,x)$ governed by a suitable diffusive prey-predator system, $u(t,x)$ and $v(t,x)$ vanish on the moving boundaries $\partial\Omega(t)$ and $\partial\Sigma(t)$, respectively. We want to understand the dynamics/variations of these two species and free boundaries. For simplicity, we assume that the interaction between these two species obeys the Lokta-Volterra law, and restrict our problem to the one dimensional case. Moreover, we think that the left boundaries of $\Omega(t)$ and $\Sigma(t)$ are fixed and coincident. So, we can take $\Omega_0=(0,g_0)$, $\Sigma_0=(0,h_0)$ with $0<h_0\le g_0$, and $\Omega(t)=(0,g(t))$, $\Sigma(t)=(0,h(t))$. Based on the {\it deduction of free boundary conditions} given in \cite{BDK, WZ12}, we have the following free boundary conditions
 $$g'(t)=-\beta u_x(t,g(t)), \ \ \ h'(t)=-\mu v_x(t,h(t)),$$
where positive constants $\beta=d_1k^{-1}_1$ and $\mu=d_2k^{-1}_2$ can be considered as the {\it moving parameters}, $d_1,d_2$ and $k_1,k_2$ are, respectively, their diffusion coefficients and {\it preferred density levels} nearing free boundaries. Under the suitable rescaling, the model we are concerned here becomes the following free boundary problem
 \bes
 \left\{\begin{array}{lll}
 u_t-du_{xx}=u(a-u-bv), &t>0,\ \ 0<x<g(t),\\[0.5mm]
 v_t-v_{xx}=v(1-v+cu),\ \ &t>0, \ \ 0<x<h(t),\\[0.5mm]
 u_x(t,0)=v_x(t,0)=0,\ \ \ &t\ge0,\\[0.5mm]
 g'(t)=-\beta u_x(t,g(t)),\ \ h'(t)=-\mu v_x(t,h(t)), \ &t\ge0,\\[0.5mm]
 u(t,x)=0 \ \ {\rm for}\ x\ge g(t), \ \ v(t,x)=0 \ \ {\rm for}\ x\ge h(t), \ &t\ge0,\\[0.5mm]
  u(0,x)=u_0(x)\ \ {\rm in} \ [0, g_0],\ \ v(0,x)=v_0(x)\ \ {\rm in} \ [0,h_0], \\[0.5mm]
 g(0)=g_0\ge h_0=h(0)>0,
 \end{array}\right.\lbl{1.1}
 \ees
where $a, b, c, d, g_0, h_0, \beta$ and $\mu$ are given positive constants. The initial functions $u_0(x),v_0(x)$ satisfy
 \bes\left\{\begin{array}{l}
 u_0\in C^2([0,g_0]),\  u'_0(0)=0,\  u_0(x)>0\ \ {\rm in}\, \
[0,g_0), \ \ u_0(x)=0\ \ {\rm in}\, \ [g_0,\yy),\\[.5mm]
v_0\in C^2([0,h_0]),\  v'_0(0)=0,\  v_0(x)>0\ \ {\rm in}\, \
[0,h_0), \ \ v_0(x)=0\ \ {\rm in}\, \ [h_0,\yy).
 \end{array}\right.\lbl{1.2}
 \ees

Because the two free boundaries may intersect each other, it seems very difficult to understand the whole dynamics of this model. We shall see that the problem (\ref{1.1}) possesses the multiplicity and complexity of spreading, vanishing and asymptotic propagation. The phenomena exhibited by these multiplicities and complexities seem closer to the reality.

Some related free boundary problems of competition-diffusion model with different free boundaries have been studied recently. In \cite{DW13}, Du and Wang discussed the following problem
  \bess
 \left\{\begin{array}{lll}
 u_t-u_{xx}=u(1-u-av), &t>0,\ \ -\infty<x<g(t),\\[1mm]
  v_t-dv_{xx}=rv(1-v-bu),\ \ &t>0, \ \ h(t)<x<\infty,\\[1mm]
 u=0\ \ {\rm for}\ x\ge g(t), \ \
 v=0 \ \ {\rm for}\ x\le h(t),\ \ \ &t>0, \\[1mm]
 g'(t)=-\beta u_x(t,g(t)),\ \ h'(t)=-\mu v_x(t,h(t)),\ \ &t>0\\[1mm]
 u=u_0(x)\ \ {\rm in} \ (-\infty, g_0], \ \
 v=v_0(x) \ \ {\rm in} \ [h_0,\infty), \ \ &t=0,\\[1mm]
 g(0)=g_0<h_0=h(0).
 \end{array}\right.\eess
In this model, the competing species $u$ and $v$ occupied habitats $(-\infty,g_0]$ and $[h_0,\infty)$ at the initial time, respectively. They will move outward along free boundaries as time increases ($u$ moves to right, $v$ moves to left). When their habitats overlap, they obey the Lotka-Volterra competition law in the common habitat. Guo \& Wu \cite{GW15}, Wang \cite{Wcom15} and Wu \cite{Wu15} studied a two-species competition-diffusion model with two free boundaries, in there the left boundary conditions and free boundary conditions are the same as that of (\ref{1.1}).

The same spreading mechanism as in \cite{DLin} has been adopted in studying some two-species competition systems or prey-predator systems. The authors of \cite{DL2, DWZ15, WJ, WZh, ZhaoW} investigated a competition model in which the invasive species exists initially in a ball and invades into the new environment, while the resident species distributes in the whole space $\mathbb{R}^N$. In \cite{WZ12}, Wang and Zhao studied a predator-prey model with double free boundaries in which the predator exists initially in a bounded interval and invades into the new environment from two sides, while the prey distributes in the whole line $\mathbb{R}$.  In \cite{GW12, WZjdde}, two competition species are assumed to spread along the same free boundary. Predator-prey models with homogeneous Dirichlet (Neumann, Robin) boundary conditions at the left side and free boundary at the right side can be found in \cite{Wjde14, Wcnsns15, WZhang, ZhangW}. For traveling wave solutions of free boundary problems, see \cite{CC12, CF, Wu13} for examples.

There have been many papers concerning the free boundary problems of single equation to describe the spreading mechanism of an invading species. Please refer to \cite{DG}-\cite{DLin}, \cite{DLou}-\cite{DMW}, \cite{GLZ, KM, KY, LL1, Wjde15, Wperiodic15} and the references therein.

This paper is organized as follows. In Section 2 we study the global existence, uniqueness, regularity and some estimates of $(u,v, g, h)$. In section 3, we first recall some fundamental results from \cite{BDK, DLin} and then give some rough
estimates, which will be used in the following two sections. Section 4 is concerned with the long time behaviors of $(u,v)$, and Section 5 deals with conditions for spreading and vanishing. In Section 6, we provide some estimates of asymptotic speeds of $g(t)$, $h(t)$ and asymptotic spreading speeds of $u(t,x)$, $v(t,x)$. Finally, in section 7 we give a brief discussion.

\section{Existence, uniqueness and estimates of global solution}
\setcounter{equation}{0} {\setlength\arraycolsep{2pt}

In this section, we first prove the following local existence and uniqueness results. Then we give some uniform estimates and show that the solution exists globally in time $t$. The main ideas of this article and literature \cite{GW15} are to straighten the free boundary and use the fixed point theorem. However, in \cite{GW15} the authors considered a map of $(g(t),h(t))$ and used the contraction mapping theorem directly. In the present paper, based on the results of single equation,  we shall use Schauder's fixed point theorem to get the existence of local solution and then prove the uniqueness.

In order to facilitate the writing, we denote
  \[\Lambda=\big\{a,b,c,d,g_0,h_0,\beta,\mu,\|u_0\|_{W_p^2},
  \|v_0\|_{W_p^2}\big\}.\]
For the given positive constants $T$, $m$ and function $f(t)$, we set
 \bess
 \Delta_T^m=[0,T]\times[0,m], \ \ D_T^f=\{0\le t\le T, \ 0\le x<f(t)\}, \ \ D_\infty^f=\{t>0, \ 0\le x\le f(t)\}.\eess

\begin{theo}\lbl{th2.1}\, For any given $(u_0,v_0)$ satisfying $(\ref{1.2})$, $\alpha\in(0,1)$ and $p>3/(1-\alpha)$, there is a $T>0$ such that the problem $(\ref{1.1})$ admits a unique solution
  $$(u,v,g,h)\in W_p^{1,2}(D^g_T)\times W_p^{1,2}(D^h_T)\times
  [C^{1+\frac \alpha 2}([0,T])]^2.$$
Moreover,
  \bes
 \|u\|_{W_p^{1,2}(D^g_T)}+\|v\|_{W_p^{1,2}(D^h_T)}+\|g\|_{C^{1+\frac \alpha
2}([0,T])}+\|h\|_{C^{1+\frac \alpha 2}([0,T])}\leq C,\lbl{2.1}
  \ees
where positive constants $C$ and $T$ depend only on $\Lambda$, $\alpha$ and $p$.
 \end{theo}

{\bf Proof}. Some techniques of this proof have been inspired by \cite[Theorem 2.1]{DW13}.  We will divide the proof into three steps.

{\it Step 1: Transformation of the problem \eqref{1.1}}. Let
 \[y=x/g(t), \ \  w(t,y)=u\big(t,g(t)y\big), \ \ z(t,y)=v\big(t,g(t)y\big), \]
then \eqref{1.1} is equivalent to
  \bes
 \left\{\begin{array}{lll}
 w_t-\frac d {g^2(t)}w_{yy}-\frac{g'(t)}{g(t)}y w_y=w(a-w-b z),\ &0<t\le T, \ \ 0<y<1,\\[1mm]
 w_y(t,0)=w(t,1)=0,\ \  &0\le t\le T, \\[0.5mm]
w(0,y)=u_0(g_0y):=w_0(y), \ &0\leq y\le 1,\\[.5mm]
 g'(t)=-\beta\frac 1{g(t)}w_y(t,1),\ \ & 0\le t\le T,\\[.5mm]
  g(0)=g_0,\end{array}\right.
  \lbl{2.3}\ees\vspace{-4mm}
 \bes\left\{\begin{array}{lll}
 z_t-\frac 1 {g^2(t)}z_{yy}-\frac{g'(t)}{g(t)}y z_y=z(1-z+cw),\ &0<t\le T, \ \ 0<y<\frac{h(t)}{g(t)},\\[1.5mm]
 z_y(t,0)=z(t,\frac{h(t)}{g(t)})=0,\ \ &0\le t\le T, \\[0.5mm]
 z(0,y)=v_0(g_0y):=z_0(y), \ &0\leq y\le h_0/g_0,\\[0.5mm]
 h'(t)=-\mu\frac 1{g(t)}z_y(t,s(t)),\ \ \ &0\le t\le T,\\[0.5mm]
 \ h(0)=h_0.
 \end{array}\right.\lbl{2.5}
 \ees

{\it Step 2: Existence of the solution $(w,z,g,h)$ to \eqref{2.3} and \eqref{2.5}}.
Denote $m=1+h_0/g_0$ and define $\hat z_0(y)=z_0(y)$ for $0\le y\leq h_0/g_0$, $\hat z_0(y)=0$ for $h_0/g_0\le y\le m$. For $0<T\ll 1$, we set
 \bess
 Z_T=\big\{z\in C(\Delta_T^m):\,z(0,y)=\hat z_0(y), \,
 \|z-\hat z_0\|_{C(\Delta_T^m)}\leq 1\big\}.
 \eess
Then $Z_T$ is a bounded and closed convex set of $C(\Delta_T^m)$. For the given $z\in Z_T$, we consider $z=z(y,t)$ as a coefficient. Since $z$ satisfies
 \[\|z\|_{C(\Delta_T^1)}\le\|z\|_{C(\Delta_T^m)}\le 1+\|\hat z\|_{C(\Delta_T^m)}\le 1+\|v_0\|_{L^\infty},\]
similarly to the arguments in the proof of \cite[Theorem 2.1]{DLin} (\cite[Theorem 2.1]{WZ12}), by using the contraction
mapping theorem we can prove that, when $0<T\ll 1$, the
problem \eqref{2.3} has a unique solution $(w(t, y), g(t))$ and $w\in W^{1,2}_p(\Delta_T^1)\hookrightarrow C^{\frac{1+\alpha} 2,\,1+\alpha}(\Delta_T^1)$, $g\in C^{1+\frac\alpha 2}([0,T])$. Moreover,
  \bes
 &g'(0)=-\beta g^{-1}_0w_y(0,1)=-\beta u'_0(g_0),&\nonumber \\[.6mm]
  &g'\ge 0, \ g\le g_0+1 \ \ {\rm in}\, \ [0,T]; \ \ \
  w>0 \ \ {\rm in}\, \ [0,T]\times[0,1),&\nonumber\\[.6mm]
& \|w\|_{W^{1,2}_p(\Delta_T^1)}+\|w\|_{C^{\frac{1+\alpha} 2,\,1+\alpha}(\Delta_T^1)}
 +\|g\|_{C^{1+\frac\alpha 2}([0,T])}\leq C_1,&
 \lbl{2.7}\ees
where $C_1$ depending only on $a, b, d$, $\beta$, $g_0$, $\alpha$, $p$, $\|u_0\|_{W^2_p}$ and $1+\|v_0\|_{L^\infty}$. Besides, let $z_i\in Z_T$ and $(w_i,g_i)\in W^{1,2}_p(\Delta_T^1)\cap C^{1+\frac\alpha 2}([0,T])$ be the unique solution of \eqref{2.3} with $z=z_i$. Similar to the following proof of the uniqueness we can get the estimate
 \bes
  \|w_1-w_2\|_{C(\Delta_T^1)}+\|g_1-g_2\|_{C^1([0,T])}
  \le \|z_1-z_2\|_{C(\Delta_T^1)}\lbl{x1}\ees
provided $0<T\ll 1$, this implies that $(w,g)$ depends continuously on $z$. We shall prove \eqref{x1} in the appendix.

For such a $(w(t,y),g(t))$, determined uniquely by the above, we put $w(t,y)$ zero extension to $[0,T]\times[1,m]$ and consider the problem \eqref{2.5}. Set $s(t)=\frac{h(t)}{g(t)}$, and
 \[\xi=y/s(t), \ \ \phi(t,\xi)=z(t,s(t)\xi), \ \ \psi(t,\xi)=w(t,s(t)\xi).\]
Then \eqref{2.5} is equivalent to the following problem
  \bes
\left\{\begin{array}{lll}
 \phi_t-\frac 1 {h^2(t)}\phi_{\xi\xi}-\frac{h'(t)}{h(t)}\xi \phi_\xi=\phi[1-\phi+c\psi(t,\xi)],\ &0<t\le T, \ \ 0<\xi<1,\\[1mm]
 \phi_\xi(t,0)=\phi(t,1)=0,\ \ &0\le t\le T, \\[0.5mm]
 \phi(0,\xi)=v_0(h_0\xi):=\phi_0(\xi), \ &0\leq \xi\le 1,\\[.5mm]
 h'(t)=-\mu\frac 1{h(t)}\phi_\xi(t,1),\ \ \ &0\le t\le T,\\[.5mm]
  h(0)=h_0.
  \end{array}\right.
  \label{2.10}\ees
Similarly to the above, we can prove that, when $0<T\ll 1$, the problem \eqref{2.10} has a unique solution $(\phi(t,\xi),h(t))$ which depends continuously on $(w,g)$, and thus continuously dependent on $z$. Let $\bar z(t,y)=\phi(t,\frac y{s(t)})$. Then $(\bar z(t,y), h(t))$ is the unique solution of \eqref{2.5}, and $(\bar z(t,y), h(t))$ is continuous with respect to $z$. Moreover, the following hold:

(i) $h'(0)=-\mu v'_0(h_0)$, $h'\ge0$, $h\le h_0+1$ and $h/g<1+h_0/g_0$ in $[0,T]$, $\bar z>0$ in $D_T^{h/g}$.

(ii) There exists a constant $C_2>0$ depending only on $\Lambda$, $\alpha$ and $p$, such that
 \bes \|\bar z\|_{W^{1,2}_p(D_T^{h/g})}+\|\bar z\|_{C^{\frac{1+\alpha} 2,\,1+\alpha}(\bar D_T^{h/g})}+\|h\|_{C^{1+\frac\alpha 2}([0,T])}\leq C_2.
 \lbl{2.12}\ees

Put $\bar z(t,\cdot)$ zero extension to $[h(t)/g(t),m]$ for each $t\in [0,T]$. Then $\bar z_y\in L^\infty(\Delta_T^m)$. In view of $h\in C^{1+\frac\alpha 2}([0,T])$, we can verify that $\bar z\in C^{\frac\alpha 2,\,\alpha}(\Delta_T^m)$ and, upon using (\ref{2.12}), that
 \bes
  \|\bar z_y\|_{L^\infty(\Delta_T^m)}+\|\bar z\|_{C^{\frac\alpha 2,\,\alpha}(\Delta_T^m)}\leq C_3.
 \lbl{2.13}\ees

Define a map
 \[{\cal G}:\,Z_T\longrightarrow C(\Delta_T^m), \ \ \ {\cal G}(z)=\bar z.\]
From the above arguments we see that ${\cal G}$ is  continuous in $Z_T$, and $z\in Z_T$ is a fixed point of ${\cal G}$ if and only if $(w,z,g,h)$ solves (\ref{2.3}) and \eqref{2.5} for $0<t\le T$, where $(w,g)$ is the solution of \eqref{2.3}, and $(z,h)$ is the solution of \eqref{2.5} with the zero extension of $w(t,y)$ to $[0,T]\times[1,m]$. Estimation (\ref{2.13}) indicates that ${\cal G}$ is compact.

Notice $\bar z(0,x)=\hat z_0(y)$. Using the mean value theorem and (\ref{2.13}) we have
 \bess
 \|\bar z-\hat z_0\|_{C(\Delta_T^m)}
 \leq\|\bar z\|_{C^{\frac\alpha 2,0}(\Delta_T^m)}T^{\frac\alpha 2}\leq C_3T^{\frac\alpha 2}.
 \eess
Therefore, if we take $0<T\ll 1$, then ${\cal G}$ maps $Z_T$ into itself. Hence, ${\cal G}$ has at least one fixed point $z\in Z_T$, i.e., (\ref{2.3}) and \eqref{2.5} has at least one solution $(w,z,g,h)$ in $[0,T]$. Moreover, from the above discussion we see that $(w,z,g,h)$ satisfies
\bess
&g,\,h\in C^{1+\frac\alpha 2}([0,T]),\ \ \ g'(t)\ge 0,\, \ h'(t)\ge 0 \ \ {\rm in} \ [0,T],& \\[.5mm]
& w\in W^{1,2}_p(\Delta_T^1)\cap C^{\frac{1+\alpha} 2,\,1+\alpha}(\Delta_T^1), \ \ z\in W^{1,2}_p(\bar D_T^{h/g})\cap C^{\frac{1+\alpha} 2,\,1+\alpha}(\bar D_T^{h/g}).&\eess

{\it Step 3: Existence and uniqueness of the solution $(u,v,g,h)$ to \eqref{1.1}}. Define
 \[u(t,x)=w(t,x/g(t)), \ \ \ v(t,x)=z(t,x/g(t)).\]
Then $(u,v,g,h)$ solves (\ref{1.1}), and $(u,v)$ satisfies
\bess
u\in W^{1,2}_p(D_T^g)\cap C^{\frac{1+\alpha} 2,\,1+\alpha}(\bar D_T^g), \ \ v\in W^{1,2}_p(D_T^h)\cap C^{\frac{1+\alpha} 2,\,1+\alpha}(\bar D_T^{h}).\eess
It is easy to see that (\ref{2.1}) holds.

In the following we prove the uniqueness. Let $(u_i,v_i,g_i,h_i)$, with $i=1,2$, be two solutions of (\ref{1.1}), which are defined for $t\in [0,T]$ and $0<T\ll 1$. We can think of that
 \bess
 &g_0\le g_i(t)\le g_0+1, \ \ \ h_0\le h_i(t)\le h_0+1 \ \ \ {\rm in} \ \, [0,T], \ \ i=1,2,&\\[.5mm]
 &\|u_i(t,x)-u_0(x)\|_{C(D_T^{g_i})}\le 1, \ \ \ \|v_i(t,x)-v_0(x)\|_{C(D_T^{h_i})}\le 1, \ \ i=1,2.&\eess

Take $k=g_0+h_0+1$. For each $t\in[0,T]$, define $u_i(t,\cdot)=0$ in $[g_i(t),k]$ and $v_i(t,\cdot)=0$ in $[h_i(t),k]$, $i=1,2$. Then $v_{ix}\in L^\infty(\Delta_T^k)$. Let
 \[w_i(t,y)=u_i\big(t,g_i(t)y\big), \ \ z_i(t,y)=v_i\big(t,g_i(t)y\big), \ \ \ 0\le t\le T, \ 0\le y\le 1.\]
Then $(w_i,g_i)$ solves \eqref{2.3} with $z=z_i$ and satisfies (\ref{2.7}). Set $W=w_1-w_2$, $G=g_1-g_2$, we have
  \bes\left\{\begin{array}{lll}
W_t-\frac d {g^2_1(t)}W_{yy}- \frac{g_1'(t)}{g_1(t)}y W_y-(a-w_1-w_2-b z_1)W\\[1mm]
\qquad
=d\kk(\frac 1{g^2_1(t)}-\frac 1{g^2_2(t)}\rr)w_{2yy}+\kk(\frac{g_1'(t)}{g_1(t)}- \frac{g_2'(t)}{g_2(t)}\rr)yw_{2y}\\[1mm]
\qquad \ \ \ -bw_2(z_1-z_2),\quad \, \ \quad 0<t\le T, \ \ 0<y<1,\\[1mm]
W_y(t,0)=W(t,1)=0,\quad\ \ 0\le t\le T,\\[1mm]
W(0,y)=0, \ \quad 0\le y\le 1,
 \end{array}\right.\lbl{2.14}\ees
and
 \bes\left\{\begin{array}{lll}
  G'(t)=-\beta\frac 1{g_1(t)}W_y(t,1)-\beta\kk(\frac 1{g_1(t)}-\frac 1{g_2(t)}\rr)w_{2y}(t,1),\quad \ 0\le t\le T,\\[2mm]
  G(0)=0.\end{array}\right.
  \lbl{2.15}\ees
Remember the facts $\|w_2\|_{W^{1,2}_p(\Delta_T^1)}\le C_1$, $g_0\le g_i(t)\le g_0+1$, $|g_i'(t)|\le C_1$ and $\|z_i\|_{C(\Delta_T^1)}\le 1+\|v_0\|_{L^\infty}$, $i=1,2$. We can apply the $L^p$ estimate to (\ref{2.14}) and use Sobolev's imbedding theorem to derive
 \bes
 \|W\|_{C^{\frac{1+\alpha} 2,\,1+\alpha}(\Delta_T^1)}
 \leq C_4 \big(\|z_1-z_2\|_{C(\Delta_T^1)}+\|G\|_{C^1([0,T])}\big).\lbl{b.10}\ees
Now we estimate $\|z_1-z_2\|_{C(\Delta_T^1)}$. For any fixed $(t,y)\in\Delta_T^1$, we have
 \bess
 |z_1(t,y)-z_2(t,y)| &\leq& |v_1(t,g_1(t)y)-v_1(t, g_2(t)y)|+|v_1(t, g_2(t)y)-v_2(t,g_2(t)y)|\\
 &\leq& \|v_{1x}\|_{L^\infty(\Delta_T^k)}\|G\|_{C([0,T])}+\|v_1-v_2\|_{C(\Delta_T^k)},
 \eess
which implies,
 \bess
 \|z_1-z_2\|_{C(\Delta_T^1)}\leq\|v_1-v_2\|_{C(\Delta_T^k)}
 +\|v_{1x}\|_{L^\infty(\Delta_T^k)}\|G\|_{C([0,T])}.\eess
Combining this with \eqref{b.10}, we get
 \bes
 \|W\|_{C^{\frac{1+\alpha} 2,\,1+\alpha}(\Delta_T^1)}
 \leq C_5\big(\|v_1-v_2\|_{C(\Delta_T^k)}+\|G\|_{C^1([0,T])}\big).\lbl{b.11}\ees
Therefore, by use of (\ref{2.15}),
  \bes\|G'\|_{C^{\frac\alpha 2}([0,T])}&\leq& \beta\|g^{-1}_1W_y\|_{C^{\frac\alpha 2,0}(\Delta_T^1)}
 +\beta\|(g_1^{-1}-g_2^{-1})w_{2y}\|_{C^{\frac\alpha 2,0}(\Delta_T^1)}\nonumber\\[.5mm]
 &\leq& C_6\big(\|v_1-v_2\|_{C(\Delta_T^k)}+\|G\|_{C^1([0,T])}\big).\lbl{b.12}\ees
Recall $W(0,y)=0$, $G(0)=G'(0)=0$. Take advantage of the mean value theorem and \eqref{b.11}, \eqref{b.12}, it follows that
 \bess
 \|W\|_{C(\Delta_T^1)}
 &\leq&T^{\frac\alpha 2}\|W\|_{C^{\frac\alpha 2,0}(\Delta_T^1)}\leq C_5T^{\frac\alpha 2}\big(\|v_1-v_2\|_{C(\Delta_T^k)}+\|G\|_{C^1([0,T])}\big),\\[1mm]
 \|G\|_{C^1([0,T])}&\leq& 2T^{\frac\alpha 2}\|G'\|_{C^{\frac\alpha 2}([0,T])}
 \leq 2C_6T^{\frac\alpha 2}\big(\|v_1-v_2\|_{C(\Delta_T^k)}+\|G\|_{C^1([0,T])}\big).\eess
Make the zero extension of $w_i(t,\cdot)$ to $[1, k/g_0]$ for each $t\in[0,T]$. The above estimates lead to
 \bes
 \|W\|_{C(\Delta_T^{k/g_0})}+\|G\|_{C([0,T])}\leq C_7T^{\frac\alpha 2}\|v_1-v_2\|_{C(\Delta_T^k)}\lbl{2.17}\ees
provided $0<T\ll 1$. Moreover, because $w_i$ satisfies (\ref{2.7}), it is easy to show that
 \bess
 w_{iy}\in L^\infty(\Delta_T^{k/g_0}), \ \ \  \|w_{iy}\|_{L^\infty(\Delta_T^{k/g_0})}\leq C_1.
 \eess

Now we estimate $\|u_1-u_2\|_{C(\Delta_T^k)}$. For any $(t,x)\in \Delta_T^k$, we have $0\le x\le k$ and
 \bess
 |u_1(t,x)-u_2(t,x)|&\leq&|w_1(t,g^{-1}_1(t)x)-w_2(t,g^{-1}_1(t)x)|
 +|w_2(t,g^{-1}_1(t)x)-w_2(t,g^{-1}_2(t)x)|\\[.8mm]
 &\le& \|w_1-w_2\|_{C(\Delta_T^{k/g_0})}
 +k\|w_{2y}\|_{L^\infty(\Delta_T^{k/g_0})}
 |g^{-1}_1(t)-g^{-1}_2(t)|\\[0.8mm]
 &\le& \|w_1-w_2\|_{C(\Delta_T^{k/g_0})}+C_8\|g_1-g_2\|_{C([0,T])},
 \eess
where $C_8=mg_0^{-2}C_1$. This implies $\|u_1-u_2\|_{C(\Delta_T^k)}\leq \|W\|_{C(\Delta_T^{k/g_0})}+C_8 \|G\|_{C([0,T])}$.
In consideration of (\ref{2.17}), it follows that
\bess
 \|u_1-u_2\|_{C(\Delta_T^k)}+\|g_1-g_2\|_{C([0,T])}\leq
C_9T^{\frac\alpha 2}\|v_1-v_2\|_{C(\Delta_T^k)}.\eess

Similarly,
 \bess
 \|v_1-v_2\|_{C(\Delta_T^k)}+\|h_1-h_2\|_{C([0,T])} &\leq&C_{10}T^{\frac\alpha 2}\|u_1-u_2\|_{C(\Delta_T^k)}.\eess
Thus, when $0<T\ll 1$, we have $(u_1,v_1,g_1,h_1)=(u_2,v_2,g_2,h_2)$. The uniqueness is obtained and the proof is finished.\ \ \ \fbox{}

\vskip 4pt To show that the local solution obtained in Theorem \ref{th2.1}
can be extended in time $t$, we need the following estimates, their proofs are similar to those of \cite[Lemma 2.1]{WZ12}, and the details will be omitted here.

\begin{lem}\label{l2.1} \, Let $T\in(0,\infty)$ and $(u,v,g,h)$ be a solution of $(\ref{1.1})$ defined in $[0,T]$. Then
  $$\begin{array}{lll}
  0<u(t,x)\leq\max\{a,\,\|u_0\|_\infty\}:=M_1, \ \ &\forall \  0\le t\leq T,\ 0\le x<g(t),\\[0.5mm]
0<v(t,x)\leq\max\{1+c M_1,\,\|v_0\|_\infty\}:=M_2, \ \ &\forall \  0\le t\leq T,\ 0\le x<h(t),\\[2mm]
 0\le g'(t)\leq 2\beta\max\kk\{M_1\sqrt{a/(2d)},\,-\min_{[0,\,g_0]}u_0'(x)\rr\},\ &\forall \ 0\le t\leq T,\\[2mm]
0\le h'(t)\leq 2\mu\max\kk\{M_2\sqrt{(1+c M_1)/2},\,-\min_{[0,\,h_0]}v_0'(x)\rr\}, \ &\forall \ 0\le t\leq T.
  \end{array}$$
\end{lem}

\begin{theo}\lbl{th2.2} \ The problem $(\ref{1.1})$ admits a unique global solution $(u,v,g,h)$ and $g'(t), h'(t)>0$,
    \bes
(u,v,g,h)\in C^{1+\frac{\alpha}2,\,2+\alpha}(D_\infty^g)\times C^{1+\frac{\alpha}2,\,2+\alpha}(D_\infty^h)\times
 [C^{1+\frac{1+\alpha}2}([0,\infty))]^2.
 \lbl{2.18}\ees

Moreover, there exists a positive constant $C$, depends only on $\Lambda$, such that
  \bes
  \|u(t,\cdot)\|_{C^1([0,\,g(t)])}\leq C, \ \ \forall \ t\ge
  1; \ \ \ \ \|g'\|_{C^{{\alpha}/2}([1,\yy))}\leq C,
  \lbl{2.19}\\[1mm]
 \|v(t,\cdot)\|_{C^1([0,\,h(t)])}\leq C, \ \ \forall \ t\ge 1; \ \ \ \ \|h'\|_{C^{{\alpha}/2}([1,\yy))}\leq C.
  \lbl{2.20}\ees
 \end{theo}

{\bf Proof}. Thanks to the estimate (\ref{2.1}) and Lemma \ref{l2.1}, we can extend the unique local solution $(u,v,g,h)$ obtained in Theorem \ref{th2.1} to a global solution and
  \bes
u\in C^{\frac{1+\alpha} 2,\,1+\alpha}(D_\infty^g),\ \ v\in C^{\frac{1+\alpha} 2,\,1+\alpha}(D_\infty^h), \ \ g,\,h\in C^{1+\frac\alpha 2}([0,\infty)),
  \lbl{2.21}\ees
see \cite{DW13, GW15} for the details. Moreover, make use of \cite[Lemma 2.6]{Lie}, it can be deduced that $g'(t)>0,\, h'(t)>0$ for $t>0$. Note that $g,\,h\in C^{1+\alpha/ 2}([0,\infty))$. It is easy to verify that, for any given $T, k>0$, $u,v\in C^{\alpha/2,\,\alpha}(\Delta_T^k)$.
Using this fact and (\ref{2.21}), we can prove (\ref{2.18}) by the similar way to that of \cite[Theorem 2.1]{Wperiodic15}. The details are omitted here.

Now we prove \eqref{2.19}. It suffices to show that
 \bes
 \|u(t,\cdot)\|_{C^1([0,\,g(t)])}\leq C, \ \ \forall \ t\ge
  1; \ \ \ \ \|g'\|_{C^{{\alpha}/2}([n+1,n+3])}\leq C', \ \ \forall \ n\geq 0.
  \lbl{2.19a}\ees
In fact, the second estimate implies that $g'(t)$ is bounded in $[1,\yy)$. This combined with \eqref{2.19a} allows us to derive $|g'(t+\sigma)-g'(t)|\le C\sigma^{\ap/2}$ for some
constant $C>0$ and all $t\ge 1$, $\sigma\ge 0$. Hence \eqref{2.19} holds.

When $g_\yy<\yy$, similarly to the arguments of \cite[Theorem 2.1]{Wperiodic15} we can obtain \eqref{2.19a}. In the following we consider the case $g_\yy=\yy$. For the integer $n\ge 0$, let $u^n(t,x)=u(n+t,x)$. Then $u^n$ satisfies
 \bes
\left\{\begin{array}{ll}
u^n_t-d u^n_{xx}-f^n(t,x)u^n=0, \ &0<t\le 3,\ 0<x<g(n+t),\\[2mm]
u^n_x(t,0)=0, \ \ u^n(t,g(n+t))=0,\ \ &0\le t\le 3,\\[2mm]
u^n(0,x)=u(n,x),&0\leq x\leq g(n),
\end{array}\right.
 \lbl{2.22}\ees
where
 \[f^n(t,x)=a-u(n+t,x)-bv(n+t,x).\]
According to Lemma \ref{l2.1}, we know that $u^n$ and $f^n$ are bounded uniformly on $n$, and $g(n+t)\leq g(n+1)+M(t-1)\leq g(n+1)+2M$ for $1\le t\le 3$, where
  \[M=2\beta\max\kk\{\sqrt{\frac a{2d}}\max\{a,\,\|u_0\|_\infty\},\,-\min_{[0,\,g_0]}u_0'(x)\rr\}.\]
As $g_\yy=\yy$, there exists an $n_0\ge 0$ such that
 \[g(n_0+1)>2M+2, \ \ g(n_0)>3.\]
In the same way as the proof of \cite[Theorem 2.1]{Wperiodic15} we can show that
\bes
  \|u(t,\cdot)\|_{C^1([0,\,g(t)])}\leq C, \ \ \forall \ 1\le t\le
 n_0+3; \ \ \ \ \|g'\|_{C^{{\alpha}/2}([n+1,n+3])}\leq C, \ \ \forall \ n\leq n_0.
  \lbl{2.23}\ees

Choose $p\gg 1$. For any integer $0\le k\le g(n+1)-3$, we can apply the interior $L^p$ estimate (cf. \cite[Theorem 7.20]{Lie}) to the problem \eqref{2.22} and derive that there exists a positive constant $C$ independent of $k$ and $n$ such that
  \[\|u^n\|_{W^{1,2}_p([1,3]\times[k,k+2])}\le C, \ \ \forall \ k,n\ge 0.\]
By the embedding theorem, $\|u^n\|_{C^{\frac{1+\alpha}2, 1+\alpha}([1,3]\times[k,k+2])}\le C$, which leads to $\|u^n\|_{C^{0, 1}([1,3]\times[k,k+2])}\le C$ for all $n\ge n_0$ and $0\le k\le g(n+1)-3$. Since these intervals $[k,k+2]$ overlap and $C$ is independent of $k$, it follows that $\|u^n\|_{C^{0, 1}([1,3]\times[0,[g(n+1)]-1])}\le C$. Therefore
 \bes
 \|u^n\|_{C^{0, 1}([1,3]\times[0,g(n+1)-2])}\le\|u^n\|_{C^{0, 1}([1,3]\times[0,[g(n+1)]-1])}\le C, \ \  \forall\ n\ge n_0,\lbl{2.24}\ees
here $[g(n+1)]$ is the integral part of $g(n+1)$. Notice that
 \[g(n+1)-2\ge g(n+t)-2M-2\ge g(n+t_0)-2M-2\ge g(n_0+1)-2M-2>0\]
for all $n\ge n_0$ and $1\le t\le 3$. Make use of the estimate \eqref{2.24} we get
  \bess
 \|u^n(t,\cdot)\|_{C^1([0,\,g(n+t)-2M-2])}\le C, \ \  \forall\ n\ge n_ 0, \ 1\le t\le 3.\eess
This leads to
  \bes\|u\|_{C^{0, 1}(E_n)}\le C, \ \  \forall\ n\ge n_ 0,\lbl{2.25}
 \ees
where
 \[E_n=\{(t,x):\,n+1\le t\le n+3,\ 0\le x\le g(t)-2M-2\}, \ \ n\ge n_0.\]
Since these rectangles $E_n$ overlap and $C$ is independent of $n$, it follows from \eqref{2.25} that
 \bes\|u(t,\,\cdot)\|_{C^1([0,\,g(t)-2M-2])}\le C, \ \ \forall\ t\ge n_0+1.\lbl{2.26}
 \ees

In the following we shall show that
 \bes
  \|u(t,\cdot)\|_{C^1([g(t)-2M-2,\,g(t)])}\leq C, \ \ \forall \ t\ge
  n_0+1; \ \ \ \ \|g'\|_{C^{{\alpha}/2}([n+1,n+3])}\leq C, \ \ \forall \ n\geq n_0.
  \lbl{2.27}\ees
Once this is done, using \eqref{2.23} and \eqref{2.26} we can derive \eqref{2.19a}.

Let $y=g(t)-x$ and $w(t,y)=u(t,g(t)-y)$. Then $w(t,y)$ satisfies
 \bess
\left\{\begin{array}{ll}
w_t-d w_{yy}+g'(t)w_y-F(t,y)w=0, \ &0<t<\infty,\ 0<y<g(t),\\[2mm]
w(t,0)=0,\ \ w_y(t,g(t))=0, \ \ &0\le t<\infty,\\[2mm]
w(0,y)=u(0,g_0-y),&0\leq y\leq g_0,
\end{array}\right.
 \eess
where $F(t,y)=a-u(t,g(t)-y)-bv(t,g(t)-y)$. Similar to the above, for the integer $n\ge n_0$, let $w^n(t,y)=w(n+t,y)$. Then $w^n$ satisfies
\bes
\left\{\begin{array}{ll}
w^n_t-d w^n_{yy}+g'(n+t)w^n_y-F^n(t,y)w^n=0, \ &0<t\le 3,\ 0<y<g(n+t),\\[2mm]
w^n(t,0)=0,\ \ w^n_y(t,g(n+t))=0, \ \ &0\le t\le 3,\\[2mm]
w^n(0,y)=w(n,y),&0\leq y\leq g(n),
\end{array}\right.\lbl{2.28}
 \ees
where $F^n(t,y)=F(n+t,y)$. It follows from Lemma \ref{l2.1} that $w^n$, $g'(n+t)$ and $F^n$ are bounded uniformly on $n$. Remember
 \[g(n+t)-(2M+2)\ge g(1+n_0)-(2M+2)>0, \ \ \forall \ n\ge n_0, \ 1\le t\le 3.\]
For $\Omega=[1,3]\times[0,2M+2]$, applying the interior $L^p$ estimate to \eqref{2.28} and embedding theorem we have that $\|w^n\|_{C^{\frac{1+\alpha}2, 1+\alpha}(\Omega)}\le C$ for all $n\ge n_0$. Hence,  $\|w\|_{C^{\frac{1+\alpha}2, 1+\alpha}(\Omega_n)}\le C$, where $\Omega_n=[n+1,n+3]\times[0,2M+2]$, $n\ge n_0$. This fact combined with
 \[g'(t)=-\mu u_x(t,g(t))=\beta w_y(t,0),\]
allows us to derive the second estimate of \eqref{2.27}.

Obviously, $\|w\|_{C^{0, 1}(\Omega_n)}\le C$. Since these rectangles $\Omega_n$ overlap and $C$ is independent of $n$, it follows that
$\|w\|_{C^{0, 1}([n_0+1,\infty)\times[0,2M+2])}\le C$. Notice that $0\le y\le 2M+2$ is equivalent to $g(t)-2M-2\le x\le g(t)$, and $u_x(t,x)=-w_y(t,y)$, the first estimate of \eqref{2.27} is followed.

The estimate \eqref{2.20} can be proved by the similar way.
\ \ \ \fbox{}

It is worth stressing that, in general, the
smoothness of the solution cannot be further promoted because of
nonlinear source terms $u(a-u-bv)$ and $v(1-v+cu)$ are only H\"{o}lder continuous in $D_\infty^g$ and $D_\infty^h$, respectively. For example, if $h(t)<g(t)$ in $[0, T]$ for some $T>0$, then $v_x(t,h(t))<0$ and $v_x(t,x)\equiv 0$ for $x>h(t)$. Therefore, $v_x(t,x)$ is not continuous at $x=h(t)$, so is $(u(a-u-bv))_x(t,x)$.

\section{Preliminaries}
\setcounter{equation}{0} {\setlength\arraycolsep{2pt}

To establish the long time behaviors of $(u,v)$ and conditions for spreading and vanishing, in this section we will state some known results.

We first consider the logistic equation with a free boundary
 \bes
 \left\{\begin{array}{lll}
 z_t-d z_{xx}=z(\theta-z), &t>0,\ \ 0<x<\rho(t),\\[.2mm]
 z_x(t,0)=0,\ \ z(t,\rho(t))=0,\ \ &t\ge 0,\\[.2mm]
 \rho'(t)=-\gamma z_x(t,\rho(t)),&t\ge 0,\\[.2mm]
 \rho(0)=\rho_0,\ \ z(0,x)=z_0(x),&0\le x\le\rho_0,
 \end{array}\right.\label{3.1}
 \ees
where $d,\,\theta,\,\gamma$ and $\rho_0$ are positive constants. Utilize the results of \cite{DLin}, the problem (\ref{3.1}) has a unique global solution and $\lim_{t\to\infty}\rho(t)=\rho_\infty$ exists. Moreover, the following facts are true:

(a)\, If $\rho_0\geq\frac{\pi}2\sqrt{d/\theta}$, then $\rho_\infty=\infty$ for all $\gamma>0$;

(b)\, If  $\rho_0<\frac{\pi}2\sqrt{d/\theta}$, then there exists a positive constant $\gamma(d,\theta,\rho_0,z_0)$ such that $\rho_\infty=\infty$ if $\gamma>\gamma(d,\theta,\rho_0,z_0)$, while $\rho_\infty<\infty$ if $\gamma\leq\gamma(d,\theta,\rho_0,z_0)$. By use of the comparison principle we can see that $\gamma(d,\theta,\rho_0,z_0)$ is decreasing in $\theta,\rho_0$ and $z_0(x)$;

(c)\, If  $\rho_\infty=\infty$, then $\lim_{t\to\infty}z(t,x)=\theta$ uniformly in the compact subset of $[0,\infty)$.

Denote
  \bes
 \beta^*=\gamma(d,a,g_0,u_0), \ \ \mu^*=\gamma(1,1,h_0,v_0),\ \ \mu_*=\gamma(1,1+ac,h_0,v_0).\lbl{3.2}\ees

Next, we consider the problem
\bes\left\{\begin{array}{l}
 d q''-k q'+q(\theta-q)=0,\ \ 0<y<\infty,\\[1mm]
q(0)=0,\ \ q'(0)=k/\nu, \ \ q(\infty)=\theta,\\[1mm]
 k\in(0,2\sqrt{\theta d}); \ \ q'(y)>0,\ \ 0<y<\infty,
\end{array}\right.\lbl{3.3}
  \ees
where $\nu, d, \theta$ and $k$ are constants.

\begin{prop}\lbl{p3.1}{\rm(\cite{BDK, DLin})} For any given $\nu, d, \theta>0$, the problem $(\ref{3.3})$ has a unique solution $(q(y),k)$. Denote $k=k(\nu,d,\theta)$. Then $k(\nu,d,\theta)$ is strictly increasing in $\nu$ and $\theta$, respectively. Moreover,
 \bes
 \lim_{\frac{\theta\nu}{ d}\to\infty}\frac{k(\nu,d,\theta)}{\sqrt{\theta d}}=2, \ \ \
 \lim_{\frac{\theta\nu}{ d}\to 0}\frac{k(\nu,d,\theta)}{\sqrt{\theta d}}\frac{d}{\theta\nu}=\frac 1{\sqrt{3}}.
 \lbl{3.4}\ees
\end{prop}

Let $w(t)$ be the unique solution of
 \[w'=w(a-w), \ \ t>0; \ \ \ w(0)=\|u_0\|_{L^\infty}.\]
Then $w(t)\to a$ as $t\to\infty$. The comparison principle leads to
 \bess
 \limsup_{t\to\infty}u(t,x)\le a\ \ {\rm uniformly\ in} \ [0,\infty).
 \eess
Similarly,
 \bess
 \limsup_{t\to\infty}v(t,x)\le 1+ac\ \ {\rm uniformly\ in} \ [0,\infty).
  \eess
Consequently, for any given $0<\ep\ll 1$, there exists $T\gg 1$ such that
 \bes\kk\{\begin{array}{ll}
  u[a-b(1+ac+\ep)-u]\le u(a-u-bv)\le u(a-u)\ \  &{\rm in}\ \ [T,\infty)\times[0,\infty),\\[.5mm]
  v(1-v)\le v(1-v+cu)\le v[1+c(a+\ep)-v]\ \ &{\rm in}\ \ [T,\infty)\times[0,\infty).
 \end{array}\rr.\lbl{3.5}\ees
In view of \cite[Theorem 4.2]{DLin} and the comparison principle, it can be deduced that
 \bes
  &&\dd\limsup_{t\to\infty}\frac{g(t)}t\le k(\beta,d,a):=\bar k_\beta, \lbl{3.6}\\[.5mm]
 &&\liminf_{t\to\infty}\frac{h(t)}t\ge k(\mu,1,1):=\ud k_\mu,
  \lbl{3.7}\ees
 and
 \bess
 &\dd \ \ \limsup_{t\to\infty}\frac{h(t)}t\le k(\mu,1,1+c(a+\ep)),&\\[.2mm]
 &\dd\liminf_{t\to\infty}\frac{g(t)}t\ge k(\beta,d,a-b(1+ac+\ep))\ \ \ {\rm if} \ \ a>b(1+ac).&
 \eess
The arbitrariness of $\ep$ yields that
 \bes
 &\dd\limsup_{t\to\infty}\frac{h(t)}t\le k(\mu,1,1+ac):=\bar k_\mu,&\lbl{3.8}\\[.5mm]
 &\dd\liminf_{t\to\infty}\frac{g(t)}t\ge k(\beta,d,a-b(1+ac)):=\ud k_\beta\ \ \ {\rm if} \ \ a>b(1+ac).&\lbl{3.9}
 \ees

\section{Long time behavior of $(u,v)$}
\setcounter{equation}{0} {\setlength\arraycolsep{2pt}

This section concerns with the limits of $(u(t,x),\,v(t,x))$ as $t\to\infty$. We first give a lemma.

\begin{lem}\label{l7.2}\, Let $d$, $C$, $\mu$ and $m_0$ be positive constants, $w\in W^{1,2}_p((0,T)\times(0,m(t)))$ for some $p>1$ and any $T>0$, and $w_x\in C([0,\infty)\times[0,m(t)])$, $m\in C^1([0,\infty))$. If $(w,m)$ satisfies
  \vspace{-.5mm}\bess\left\{\begin{array}{lll}
 w_t-dv_{xx}\geq -C w, &t>0,\ \ 0<x<m(t),\\[1.5mm]
 w\ge 0,\ \ \ &t>0, \ \ x=0,\\[1.5mm]
 w=0,\ \ m'(t)\geq-\mu w_x, \ &t>0,\ \ x=m(t),\\[1.5mm]
 w(0,x)=w_0(x)\ge,\,\not\equiv 0,\ \ &x\in (0,m_0),\\[1.5mm]
 m(0)=m_0,
 \end{array}\right.\vspace{-.5mm}\lbl{7z.25}\eess
and
  \vspace{-.5mm}
  \bess
  &\dd\lim_{t\to\infty} m(t)=m_\infty<\infty, \ \ \ \lim_{t\to\infty} m'(t)=0,&
 \label{7.22}\\[1mm]
& \|w(t,\cdot)\|_{C^1[0,\,m(t)]}\leq M, \ \ \forall \ t>1&\label{7.23}
 \eess
for some constant $M>0$.  Then
  \vspace{-.5mm}\bess
 \lim_{t\to\infty}\,\max_{0\leq x\leq m(t)}w(t,x)=0.
 \lbl{zz7.25}\vspace{-1mm}\eess
\end{lem}

\noindent{\it Proof}.\, Firstly, the maximum principle gives $w(t,x)>0$ for $t>0$ and $0<x<m(t)$. Follow the proof of \cite[Theorem 2.2]{WZjdde} word by word we can prove this lemma and the details are omitted. \ \ \ $\Box$

\begin{theo}\lbl{th4.1} If $g_\infty<\infty\ (h_\infty<\infty)$, then
 \bess
 \lim_{t\to\infty}\max_{0\leq x\leq g(t)}u(t,x)=0 \ \ \kk(\lim_{t\to\infty}\max_{0\leq x\leq h(t)}v(t,x)=0\rr).\eess
\end{theo}

{\bf Proof}. Notice Lemma \ref{l2.1} and Theorem \ref{th2.2}, the conclusion can be deduced by Lemma \ref{l7.2} directly. \ \ \ \fbox{}

When $g_\infty<\infty$ (resp., $h_\infty<\infty$), we say that the species $u$ (resp., $v$) vanishes eventually. When $g_\infty=\infty$ (resp., $h_\infty=\infty$), we say that the species $u$ (resp., $v$) spreads successfully.

\begin{theo}\lbl{th4.2}
 {\rm(i)} If $g_\infty<\infty$ and $h_\infty=\infty$, then
 \bess
 \lim_{t\to\infty}v(t,x)=1\ \ {\rm uniformly\ in\ the\ compact\ subset\ of } \ [0,\infty).\eess

  {\rm(ii)} If $h_\infty<\infty$ and $g_\infty=\infty$, then
 \bess
 \lim_{t\to\infty}u(t,x)=a\ \ {\rm uniformly\ in\ the\ compact\ subset\ of } \ [0,\infty).\eess
 \end{theo}

{\bf Proof}. We only prove (i) as (ii) can be proved by the similar way. Firstly, utilize  the comparison principle and conclusions about the logistic equation, it is easy to get
 \bes
 \liminf_{t\to\infty}v(t,x)\ge 1\ \ {\rm uniformly\ in\ the\ compact\ subset\ of } \ [0,\infty).
 \lbl{4.1}\ees
For any given $0<\ep\ll 1$. Note that $\dd\lim_{t\to\infty}\max_{0\leq x\leq g(t)}u(t,x)=0$ (Theorem \ref{th4.1}) and $u(t,x)=0$ for $x>g(t)$. There exists $T\gg 1$ such that $u(t,x)<\ep$ for all $t\ge T$ and $x\in[0,\infty)$. Thus, $v$ satisfies
 \bess
 \left\{\begin{array}{lll}
 v_t-v_{xx}\le v(1+c\ep-v), &t\ge T,\ \ 0<x<h(t),\\[0.5mm]
 v_x(t,0)=0,\ \ v(t,h(t))=0,\ &t\ge T.
 \end{array}\right.
 \eess
This implies $\dd\limsup_{t\to\infty}v(t,x)\le 1+c\ep$ uniformly in $[0,\infty)$. And so,
 \[\limsup_{t\to\infty}v(t,x)\le 1\ \ {\rm uniformly\ in} \ [0,\infty).\]
Combining (\ref{4.1}), the desired result is obtained immediately. \ \ \ \fbox{}

When $u$ (resp., $v$) vanishes eventually and $v$ (resp., $u$) spreads successfully, our model formally reduces to the single species model. In such a case, the speed of $h(t)$ (resp., $g(t)$) is the same as the one given in \cite{DLin}, the sharp estimates of $v(t,x)$ and $h(t)$ ($u(t,x)$ and $g(t)$) as those investigated by \cite{DMZ, ZhaoW15}.

In the following we deal with the case that both two species spread successfully. We first give a local result.

\begin{theo}\lbl{th4.3} Assume that $g_\infty=h_\infty=\infty$.

{\rm(i)} For the weak predation case $b<\min\{a,\,1/c\}$, we denote
 \[A=\frac{a-b}{1+bc}, \ \ \ B=\frac{1+ac}{1+bc}.\]
Then
 \bess
 \lim_{t\to\infty}u(t,x)=A, \ \ \ \lim_{t\to\infty}v(t,x)=B
  \eess
uniformly in any compact subset of $[0,\infty)$;

{\rm(ii)} For the strong predation case $b\geq a$, we have
 \[\lim_{t\to\infty}v(t,x)=1,\ \ \ \ \lim_{t\to\infty}u(t,x)=0\]
uniformly in any compact subset of $[0,\infty)$.
 \end{theo}

Applying Propositions 2.1-2.3 of \cite{WZjdde}, we can prove Theorem \ref{th4.3} by the similar arguments to those of \cite[Theorems 4.3 and 4.4]{WZ12}. The details are omitted here.

Now, we are going to study the more accurate limits of $(u,v)$ as $t\to\infty$. Let $\ud k_\mu$ and $\ud k_\beta$ be given by (\ref{3.7}) and (\ref{3.9}), respectively.

\begin{theo}\lbl{th4.4} Suppose $g_\infty=h_\infty=\infty$. For the weak predation case $b<\min\{a,\,1/c\}$, if we further assume $a>b(1+ac)$, then for each $0<k_0<\min\{\ud k_\beta,\,\ud k_\mu\}$, there hold:
  \bes
 \lim_{t\to\infty}\max_{[0,\,k_0t]}\kk|u(t,\cdot)-A\rr|=0, \ \ \lim_{t\to\infty}\max_{[0,\,k_0t]}\kk|v(t,\cdot)-B\rr|=0.
 \lbl{4d.2}\ees
\end{theo}

{\bf Proof}. Some ideas in this proof are inspired by \cite[Theorem 7]{Wcom15}. To facilitate writing, for $\tau\ge 0$, we introduce the following free boundary problem
\bes
 \left\{\begin{array}{lll}
 z_t-D z_{xx}=z(\theta-z), &t>\tau,\ \ 0<x<s(t),\\[.2mm]
 z_x(t,0)=0,\ \ z(t,s(t))=0,\ \ &t\ge \tau,\\[.2mm]
 s'(t)=-\nu z_x(t,s(t)),&t\ge \tau,\\[.2mm]
 s(\tau)=s_0,\ \ z(\tau,x)=z_0(x),&0\le x\le s_0,
 \end{array}\right.\label{4.2}
 \ees
and set $\Gamma=(\tau,D,\theta,\nu,s_0)$, where $D,\,\theta,\,\nu$ and $s_0$ are positive constants. For any given constant $T\ge 0$ and function $f(t)$, we define
 \[\Omega_T^f=\{(t,x):\, t\ge T, \ 0\le x\le f(t)\}.\]

Recall $0<k_0<\min\{\ud{k}_\beta,\,\ud{k}_\mu\}$. Take advantage of (\ref{3.7}) and  (\ref{3.9}), there exist $0<\sigma_0\ll 1$ and $ t_\sigma\gg 1$ such that
 \bess
 &k_\sigma:=k_0+\sigma<\min\{\ud{k}_\beta,\,\ud{k}_\mu\}, \ \ \forall\ 0<\sigma\le\sigma_0,&\\[.5mm]
 & g(t)>k_\sigma t, \ \ \ h(t)>k_\sigma t,\ \ \ \forall\ t\ge t_\sigma, \ 0<\sigma\le\sigma_0.&
 \eess

The following proof will be divided into five steps. The method used here is the cross-iteration scheme. In order to construct iteration sequences, in the first four steps, we will prove that, for any fixed $0<\sigma<\sigma_0/5$,
  \bes
  &\dd\liminf_{t\to\infty}\min_{[0,\,k_{5\sigma} t]}v(t,\cdot)\ge 1:=\ud v_1,& \lbl{4.a1}\\[.5mm]
 &\dd\limsup_{t\to\infty}\max_{[0,\,k_{4\sigma}t]}u(t,\cdot)\le a-b:=
 \bar u_1,&\lbl{4.4}\\[.5mm]
 &\dd\limsup_{t\to\infty}\max_{[0,\,k_{3\sigma}t]}v(t,\cdot)\le 1+c\bar u_1:=
 \bar v_1,&\lbl{4.a2}\ees
and
\bes\kk\{\begin{array}{ll}
 \dd\liminf_{t\to\infty}\min_{[0,\,k_{2\sigma}t]}u(t,\cdot)\ge a-b\bar v_1:=\ud{u}_1,\\[.5mm]
 \dd\liminf_{t\to\infty}\min_{[0,\,k_{\sigma}t]}v(t,\cdot)\ge 1+c\ud{u}_1:=\ud{v}_2,
 \end{array}\rr.\lbl{4.a3}\ees
respectively. In the last step, we will construct four sequences $\{\bar u_i\}$, $\{\bar v_i\}$, $\{\ud u_i\}$ and $\{\ud v_i\}$, and derive the desired conclusion.

{\it Step 1}: As $h_\infty=\infty$, we can find a $ t_1\gg 1$ so that $h(t_1)>\pi/2$.
Let $(z_1,s_1)$ be the unique solution of (\ref{4.2}) with  $\Gamma=(t_1,1,1,\mu,h(t_1))$ and $z_0(x)=v(t_1,x)$.
Then $h(t)\ge s_1(t)$, $v(t,x)\ge z_1(t,x)$ in $\Omega_{t_1}^{s_1}$ by the comparison principle. And $s_1(\infty)=\infty$ since $s_1(t_1)>\pi/2$. Make use of \cite[Theorem 3.1]{ZhaoW15} (see also \cite{DMZ, GLZ}), we get
 \bes
 \lim_{t\to\infty}(s_1(t)-k_1 t)=\varsigma_1\in\mathbb{R}, \ \ \lim_{t\to\infty}\|z_1(t,x)- q_1(k_1 t+\varsigma_1-x)\|_{L^\infty([0,s_1(t)])}=0,\lbl{4.3}\ees
where $(q_1(y), k_1)$ is the unique solution of (\ref{3.3}) with $(\nu,d,\theta)=(\mu,1, 1)$, i.e., $k_1=k(\mu,1, 1)=\ud k_\mu$. Note $0<k_{5\sigma}<\ud{k}_\mu=k_1$, it is easy to see that $s_1(t)-k_{5\sigma} t\to\infty$ and $\min_{[0,\,k_{5\sigma} t]}(k_1 t+\varsigma_1-x)\to\infty$ as $t\to\infty$. Owing to $q_1(y)\nearrow 1$ as $y\nearrow\infty$, we have
$\min_{x\in[0,\,k_{5\sigma} t]}q_1(k_1 t+\varsigma_1-x)\to 1$ as $t\to\infty$. It then follows, upon using (\ref{4.3}), that $\min_{[0,\,k_{5\sigma} t]}z_1(t,\cdot)\to 1$ as $t\to\infty$. Thus, (\ref{4.a1}) holds because of $v\ge z_1$ in $\Omega_{t_1}^{s_1}$.

For any given $0<\ep\ll 1$, there exists $ t_2\gg 1$ such that
 \[g(t)>k_{5\sigma} t, \ \ v(t,x)\ge 1-\ep:=\lambda_\ep,\ \ \ \forall \ t\ge t_2, \ 0\le x\le k_{5\sigma} t.\]

{\it Step 2}: Obviously, $a-b\lambda_\ep>a-b>(a-b)/(1+bc)$. Make use of Theorem \ref{th4.3}(i), it follows that $u(t,0)\le a-b\lambda_\ep$ in $[t_3,\infty)$ for some $ t_3> t_2$. Thus, $u$ satisfies
   \bess
 \left\{\begin{array}{lll}
  u_t-du_{xx}\le u(a-b\lambda_\ep-u),\ \ &t\ge t_3, \ \ 0<x<k_{5\sigma} t,\\[0.5mm]
 u\le M_1, \ &t\ge t_3, \ 0<x\le k_{5\sigma} t,\\[.5mm]
 u(t,0)\le a-b\lambda_\ep,
 \end{array}\right.
 \eess
where $M_1$ is given by Lemma \ref{l2.1}. We will show that
 \bes
 \limsup_{t\to\infty}\max_{[0,\,k_{4\sigma}t]}u(t,\cdot)\le a-b\lambda_\ep.
 \lbl{4.a5}\ees
Once this is done, (\ref{4.4}) is obtained immediately because $\ep>0$ is arbitrary. To prove (\ref{4.a5}), we choose $0<\delta\ll 1$ and define
 \[f(t,x)=a-b\lambda_\ep+M_1{\rm e}^{\delta k_{5\sigma} t_3}{\rm e}^{\delta(x-k_{5\sigma} t)}, \ \
 t\geq t_3, \ \ 0\le x\le k_{5\sigma}t.\]
Evidently,
 \[\max_{[0,\,k_{4\sigma} t]}f(t,\cdot)\le a-b\lambda_\ep+M_1{\rm e}^{\delta k_{5\sigma}t_3}{\rm e}^{-\delta\sigma t}\to a-b\lambda_\ep\]
as $t\to\infty$, and
 \bess\kk\{\begin{array}{ll}
 f(t,0)>a-b\lambda_\ep, \ f(t,k_{5\sigma}t)> M_1, \ \ &t\ge t_3,\\[.5mm]
  f(t_3,x)> M_1, \ \ &0\le x\le k_{5\sigma} t_3.
  \end{array}\rr.\eess
It is easy to verify that, when $\delta(k_{5\sigma}+d\delta)\le a-b$,
  \[f_t-df_{xx}\ge f(a-b\lambda_\ep-f).\]
The comparison principle gives $u(t,x)\le f(t,x)$ for all $t\ge t_3$ and $0\le x\le k_{5\sigma}t$. Thus we have (\ref{4.a5}) and then obtain (\ref{4.4}).

There exists $ t_4> t_3$ such that
 \[h(t)>k_{4\sigma} t, \ \ u(t,x)\le \bar u_1+\ep:=\bar u_1^\ep<1, \ \ \ \forall \ t\ge t_4, \ 0\le x\le k_{4\sigma} t.\]

{\it Step 3}: The condition $a>b$ implies $1+c\bar u_1^\ep>(1+ac)/(1+bc)$. Similarly to Step 2, by use of Theorem \ref{th4.3}(i), there exists $t_5>t_4$ such that $v(t,0)\le 1+c\bar u_1^\ep$ in $[t_5,\infty)$. Thus, $v$ satisfies
   \bess
 \left\{\begin{array}{lll}
  v_t-v_{xx}\le v(1+c\bar u_1^\ep-v),\ \ &t\ge t_5, \ \ 0<x<k_{4\sigma} t,\\[0.5mm]
 v(t,0)\le 1+c\bar u_1^\ep, \ \ v\le M_2, \ &t\ge t_5, \ 0<x\le k_{4\sigma} t,
 \end{array}\right.
 \eess
where $M_2$ is given by Lemma \ref{l2.1}. In the same way as Step 2, it can be proved that
 \[\limsup_{t\to\infty}\max_{[0,\,k_{\bar\sigma}t]}v(t,\cdot)\le 1+c\bar u_1,\]
where $\bar\sigma=7\sigma/2$. So, (\ref{4.a2}) holds.

Take $t_6>t_5$ such that
 \bes
 g(t)>k_{\bar\sigma} t, \ \ v(t,x)\le \bar v_1+\ep:=\bar v_1^\ep,\ \ \ \forall \ t\ge t_6, \ 0\le x\le k_{\bar\sigma} t.\lbl{4.a9}\ees

{\it Step 4}: It is easy to see that $a-b\bar v_1^\ep>a-b(1+ac)$ since $0<\ep\ll 1$. So,
 $$k(\beta,d,a-b\bar v_1^\ep)>k(\beta,d,a-b(1+ac))=\ud k_\beta>k_{3\sigma}.$$
Owing to $k(\beta,d,a-b\bar v_1^\ep)\to 0$ as $\beta\to 0$ (cf. (\ref{3.4})), we can take $0<\beta^*<\beta$ so that $k(\beta^*,d,a-b\bar v_1^\ep)=k_{3\sigma}$. In this way, we get a function $q(y)$, where $(q(y), k_{3\sigma})$ is the unique solution of (\ref{3.3}) with $(\nu,d,\theta)=(\beta^*,d,a-b\bar v_1^\ep)$. Because of $g(t)>k_{3\sigma} t$ for all $t\ge t_6$, we can find a function $\tilde u\in C^2([0,k_{3\sigma} t_6])$ satisfying $\tilde u'(0)=\tilde u(k_{3\sigma} t_6)=0$,  $\tilde u(x)>0$ in} $[0,k_{3\sigma} t_6)$ and
   \[\tilde u(x)\le u(t_6,x), \ \ \forall \ x\in [0,k_{3\sigma} t_6].\]
Let $(z_2,s_2)$ be the unique solution of (\ref{4.2}) with $\Gamma=(t_6,d,a-b\bar v_1^\ep,\beta^*,k_{3\sigma} t_6)$ and $z_0(x)=\tilde u(x)$.
Then, using \cite[Theorem 3.1]{ZhaoW15}, we have
 \bes
 s_2(t)-k_{3\sigma} t\to \varsigma_2\in\mathbb{R}, \ \ \|z_2(t,x)- q(k_{3\sigma} t+\varsigma_2-x)\|_{L^\infty([0,s_2(t)])}\to 0\lbl{4.5}\ees
as $t\to\infty$. As $\bar\sigma>3\sigma$, in consideration of (\ref{4.a9}) and the first limit of (\ref{4.5}), we can think of $g(t)>s_2(t)$, $v(t,x)\le \bar v_1^\ep$ for all $t\ge t_6$ and $0\le x\le s_2(t)$. As a consequence, $u$ satisfies
  \[u_t-du_{xx}\ge u(a-b\bar v_1^\ep-u) \ \ \ {\rm in}\ \ \Omega_{t_6}^{s_2}.\]
Note that $z_{2x}(t,0)=u_x(t,0)=0$, $z_2(t,s_2(t))=0<u(t,s_2(t))$ in $[t_6,\infty)$ and $z_2(t_6,x)=\tilde u(x)\le u(t_6,x)$ in $[0,k_{3\sigma} t_6]$, it is deduced that $u\ge z_2$ in $\Omega_{t_6}^{s_2}$ by the comparison principle.

Since $q(y)\nearrow a-b\bar v_1^\ep$ as $y\nearrow\infty$, we see that $\lim_{t\to\infty}\min_{x\in[0,\,k_{2\sigma}t]}q(k_{3\sigma} t+\varsigma_2-x)=a-b\bar v_1^\ep$. Apply (\ref{4.5}) once again, it follows that $s_2(t)-k_{2\sigma} t\to\infty$ and
$\min_{[0,\,k_{2\sigma}t]}z_2(t,\cdot)\to a-b\bar v_1^\ep$ as $t\to\infty$. This gives the first inequality of (\ref{4.a3})
because $u\ge z_2$ in $\Omega_{t_6}^{s_2}$ and $\ep>0$ is arbitrary.

Similarly, we can prove the second inequality of (\ref{4.a3}).

{\it Step 5}: Five positive constants $\ud v_1$, $\bar u_1$, $\bar v_1$, $\ud u_1$ and $\ud v_2$ have been obtained. Now we define
 \[\bar u_i=a-b\underline{v}_i, \ \
\bar v_i=1+c\bar u_i,\ \
 \underline{u}_i=a-b\bar{v}_i, \ \
\underline{v}_{i+1}=1+c\underline{u}_i,\ \ i=2,3,\cdots.\]
Then (cf. the proof of \cite[Theorem 4.3]{WZ12})
 \[\lim_{i\to\infty}\bar u_i=\lim_{i\to\infty}\underline{u}_i
 =A, \ \ \lim_{i\to\infty}\bar v_i=\lim_{i\to\infty}\underline{v}_i
 =B.\]
Repeating the above process we can show that
 \bess
  &\dd\ud u_i\le\liminf_{t\to\infty}\min_{[0,\,k_0 t]}u(t,\cdot), \ \ \limsup_{t\to\infty}\max_{[0,\,k_0 t]}u(t,\cdot)\le\bar u_i, \ \ \forall\ i\geq 1.&\\[1mm]
  &\dd\ud{v}_i\le\liminf_{t\to\infty}\min_{[0,\,k_0 t]}v(t,\cdot), \ \ \limsup_{t\to\infty}\max_{[0,\,k_0 t]}v(t,\cdot)\le\bar v_i, \ \ \forall\ i\geq 1.&
 \eess
The proof is finished.  \ \ \ \fbox{}

\section{Conditions for spreading and vanishing}
\setcounter{equation}{0} {\setlength\arraycolsep{2pt}

In this section we will give some conditions to identify the spreading and vanishing of $u$, $v$. Throughout this section, the positive constants $\beta^*,\,\mu^*$ and $ \mu_*$ are given by (\ref{3.2}).

\begin{theo}\lbl{th5.1}\, {\rm(i)} If $g_0<\frac{\pi}2\sqrt{d/a}$ and $\beta\leq\beta^*$, then $g_\infty<\infty$;

 {\rm(ii)} If either $h_0\ge\pi/2$, or $h_0<\pi/2$ and  $\mu>\mu^*$, then $h_\infty=\infty$;

 {\rm(iii)} If $u_0(x)\le a$, $h_0<\frac\pi 2\sqrt{1/(1+ac)}$ and $\mu<\mu_*$, then $h_\infty<\infty$.
\end{theo}

This theorem can be proved directly by the comparison principle. We omit the details.

\begin{theo}\lbl{th5.2}\,
Assume that $u_0(x)\le a$, $h_0<\frac\pi 2\sqrt{1/(1+ac)}$ and $\mu<\mu_*$. If either $g_0> \frac{\pi}2\sqrt{d/a}$, or  $g_0<\frac{\pi}2\sqrt{d/a}$ and $\beta>\beta^*$, then $g_\infty=\infty$.
\end{theo}

{\bf Proof}. First, there exists $0<\ep\ll 1$ such that either

 (a) $g_0> \frac{\pi}2\sqrt{d/(a-b\ep)}$, or

 (b) $g_0<\frac{\pi}2\sqrt{d/(a-b\ep)}$ and $\beta>\beta_\ep:=\gamma(d,a-b\ep,g_0,u_0)$.\\
As $u_0(x)\le a$, $h_0<\frac\pi 2\sqrt{1/(1+ac)}$ and $\mu<\mu_*$, using Theorem \ref{th5.1}(iii) and Theorem \ref{th4.1}, successively, we have $h_\infty<\infty$ and $\lim_{t\to\infty}\max_{0\leq x\leq h(t)}v(t,x)=0$. Note that $v(t,x)=0$ for $x\ge h(t)$. We can find a $T\gg 1$ so that   $0\le v(t,x)\le\ep$ for all $t\ge T$ and $x\ge 0$. Hence, $(u,g)$ satisfies
 \bess
 \left\{\begin{array}{lll}
 u_t-du_{xx}\ge u(a-b\ep-u), &t\ge T,\ \ 0<x<g(t),\\[0.5mm]
 u_x(t,0)=u(t,g(t))=0,\ \ \ &t\ge T,\\[0.5mm]
 g'(t)=-\beta u_x(t,g(t)), &t\ge T.
 \end{array}\right.
 \eess
Consequently, $g_\infty=\infty$ because at least one of conditions (a) and (b) holds. \ \ \ \fbox{}

The conclusions of Theorem \ref{th5.2} show that, if one of $g_0$ and $\beta$ (the initial habitat and moving parameter of the prey) is ``suitably large", both $h_0$ and $\mu$ (the initial habitat and moving parameter of the predator) are ``suitably small", the prey will spread successfully, while the predator will vanishes eventually.  However, the predator always able to successfully spread if either $h_0\ge\pi/2$, or $h_0<\pi/2$ and  $\mu>\mu^*$. A natural question arises: does the prey always die out eventually if the predator spreads successfully? Intuitively, if the predator spreads faster enough than the prey, the prey would have no chance to survive eventually even its initial population
and initial habitat size are large.

In the following, we will give two results to answer the above question. The first one indicates that if the predator spreads slowly and the prey's initial habitat is much larger than that of the predator, the prey will spread successfully and its territory  always cover that of the predator no matter whether the latter successful spread. The second one illustrates that for the strong predation case $b>a$, if the prey spreads slowly and the predator spreads quickly, the prey will vanish eventually and the predator will spread successfully.

In view of Lemma \ref{l2.1}, we have
 \[0<h(t)\le K\mu t+h_0, \ \ \forall \ t>0,\]
where
\bess
 &K=2\max\big\{M_2\sqrt{(1+c M_1)/2},\,-\min_{[0,\,h_0]}v_0'(x)\big\},& \\[1mm] & M_2=\max\{1+c M_1,\,\|v_0\|_\infty\}, \ \ M_1=\max\{a,\,\|u_0\|_\infty\}.&
 \eess

\begin{theo}\lbl{th5.4} Let $d, a, b, c$ and $\beta$ be fixed. Then there exists $0<\bar\mu<\sqrt{2da}/K$ such that, when
  \[0<\mu<\bar\mu, \ \ g_0-h_0>\frac{2d\pi}{\sqrt{2da-K^2\mu^2}}:=L(\mu),\]
we have $g(t)\ge K\mu t+h_0+L(\mu)$ for all $t\ge 0$, which leads to $g(t)> h(t)$ for all $t\ge 0$ and $g(t)\to\infty$ as $t\to\infty$.

Moreover, if $h_0\ge\pi/2$, we also have $h_\infty=\infty$ for all $\mu>0$.
 \end{theo}

{\bf Proof}. This proof is similar to that of \cite[Theorem 6]{Wcom15}. For the completeness and convenience to reader, we shall give the details. Denote $\sigma=K\mu$. For these $t$ satisfying $g(t)>\sigma t+h_0$, we define
 \[y=x-\sigma t-h_0, \ \ \vp(t,y)=u(t,x), \ \ \eta(t)=g(t)-\sigma t-h_0.\]
Then $\vp(t,y)>0$ for $t\ge 0$ and $0\le y<\eta(t)$. Note that  $v(t,x)=0$ for $x\ge h(t)$, and $y\ge 0$ implies $x\ge h(t)$, it follows that $\vp$ satisfies
 \bess\kk\{\begin{array}{ll}
 \vp_t-d \vp_{yy}-\sigma \vp_y=\vp(a-\vp), \ \ &t>0, \ 0<y<\eta(t),\\[.5mm]
 \vp(t,0)=u(t,\sigma t+h_0), \ \ \vp(t,\eta(t))=0, \ \ &t\ge 0,\\[.5mm]
\vp(0,y)=u_0(y+h_0), \ \ \ &0\le y\le g_0-h_0.
 \end{array}\rr.\eess
Let $\lambda$ be the principal eigenvalue of
 \bes\left\{\begin{array}{ll}
 -d\phi''-\sigma\phi'-a\phi=\lambda\phi, \ \ 0<x<L,\\[.5mm]
 \phi(0)=0=\phi(L).
 \end{array}\right.\lbl{5.2}\ees
Then the following relation (between $\lambda$ and $L$) holds:
 \[\frac{\pi}L=\frac{\sqrt{4d(a+\lambda)-\sigma^2}}{2d}.\]
Choose $\lambda=-a/2$ and define
 \[L_\sigma=\frac{2d\pi}{\sqrt{2da-\sigma^2}}, \ \ \ \phi(y)={\rm e}^{-\frac\sigma{2d}y}\sin\frac{\pi}{L_\sigma}y \ \ \ {\rm with} \ \ 0<\sigma<\sqrt{2da}.\]
Then $(L_\sigma,\phi)$ satisfies (\ref{5.2}) with $\lambda=-a/2$ and $L=L_\sigma$.
Assume $g_0-h_0>L_\sigma$ and set
 \[\delta_\sigma=\min\kk\{\inf_{(0,\,L_\sigma)}\frac{\vp(0,y)}{\phi(y)}, \ \frac a2\inf_{(0,\,L_\sigma)}\frac 1{\phi(y)}\rr\}, \ \ \psi(y)=\delta_\sigma\phi(y).\]
Then $0<\delta_\sigma<\infty$. It is easy to verify that $\psi(y)\le \vp(0,y)$ in $[0,L_\sigma]$ and $\psi(y)$ satisfies
   \bess\left\{\begin{array}{ll}
 -d\psi''-\sigma\psi'\le \psi(a-\psi), \ \ 0<x<L_\sigma,\\[1mm]
 \psi(0)=0=\psi(L_\sigma).
 \end{array}\right.\eess
Take a maximal $\bar\sigma\in(0,\sqrt{2da})$ so that
  \bes\sigma<\beta\delta_\sigma\frac{\pi}{L_\sigma}\exp\kk(-\frac{\sigma L_\sigma}{2d}\rr), \ \ \forall \ \sigma\in(0,\bar\sigma).\lbl{5.3}\ees

For any given $\sigma\in(0,\bar\sigma)$, we claim that $\eta(t)>L_\sigma$ for all $t\ge 0$, which implies
 $$g(t)\ge \sigma t+h_0+L_\sigma\to\infty.$$
In fact, note $\eta(0)=g_0-h_0>L_\sigma$, if our claim is not true, then we can find a $t_0>0$ such that $\eta(t)>L_\sigma$ for all $0\le t<t_0$ and $\eta(t_0)=L_\sigma$. Therefore, $\eta'(t_0)\le 0$, i.e, $g'(t_0)\le \sigma$. On the other hand, by the comparison principle, we have  $\vp(t,y)\ge\psi(y)$ in $[0, t_0]\times[0,L_\sigma]$. Particularly, $\vp(t_0,y)\ge\psi(y)$ in $[0,L_\sigma]$. Due to $\vp(t_0,L_\sigma)=\vp(t_0,\eta(t_0))=0=\psi(L_\sigma)$, it derives that
 $$\vp_y(t_0,L_\sigma)\le\psi'(L_\sigma)=-\delta_\sigma\frac{\pi}{L_\sigma}
 \exp\kk(-\frac{\sigma L_\sigma}{2d}\rr).$$
It follows, upon using $u_x(t_0,g(t_0))=\vp_y(t_0,\eta(t_0))$, that
 \[\sigma\ge g'(t_0)=-\beta \vp_y(t_0,\eta(t_0))=-\beta \vp_y(t_0,L_\sigma)\ge\beta\delta_\sigma\frac{\pi}{L_\sigma}
 \exp\kk(-\frac{\sigma L_\sigma}{2d}\rr).\]
It is in contradiction with (\ref{5.3}).

Take $\bar\mu=\bar\sigma/K$, $L(\mu)=L_\sigma$. Then $0<\mu<\bar\mu$ is equivalent to $0<\sigma<\bar\sigma$, and $g_0-h_0>L(\mu)$ is equivalent to $g_0-h_0>L_\sigma$.

At last, when $h_0\ge\pi/2$, we have $h_\infty=\infty$ for any $\mu>0$ by Theorem \ref{th5.1}(ii). The proof is complete. \ \ \fbox{}

In consideration of (\ref{3.4}), it is easy to see that
 \[\lim_{\beta\to 0}k(\beta,d,a)=0, \ \ \lim_{\mu\to\infty}k(\mu,1,1)=2.\]
By the monotonicity of $k(\nu,d,\theta)$ in $\nu$, there exist $\bar\beta,\,\bar\mu>0$ such that $k(\beta,d,a)<k(\mu,1,1)$ for all  $0<\beta\le\bar\beta$ and $\mu\ge\bar\mu$. Therefore, $(0,\bar\beta]\times[\bar\mu,\infty)\subset{\cal A}$, where
 \[{\cal A}=\{(\beta,\mu):\,\beta,\mu>0,\, k(\beta,d,a)<k(\mu,1,1)\}.\]

\begin{theo}\lbl{th5.3}\, Assume that $(\beta,\mu)\in{\cal A}$. If $b>a$ and $h_\infty=\infty$, then $g_\infty<\infty$.
\end{theo}

{\bf Proof}. Firstly, because of $b>a$, there exists $0<\ep\ll 1$ such that $a<b(1-\ep)$.

There exists $t_1\gg 1$ such that $h(t_1)>\pi/2$. Let $(z_1,s_1)$ be the unique solution of (\ref{4.2}) with  $\Gamma=(t_1,1,1,\mu,h(t_1))$ and $z_0(x)=v(t_1,x)$.
Then  $s_1(\infty)=\infty$, $h(t)\ge s_1(t)$, $v(t,x)\ge z_1(t,x)$ in $\Omega_{t_1}^{s_1}$. Moreover, make use of \cite[Theorem 3.1]{ZhaoW15} (see also \cite{DMZ, GLZ}) we have that, as $t\to\infty$,
 \bes
 s_1(t)-k_1 t\to \varsigma_1\in\mathbb{R}, \ \ \|z_1(t,x)- q_1(k_1 t+\varsigma_1-x)\|_{L^\infty([0,s_1(t)])}\to 0,\lbl{5.1}\ees
where $(q_1(y), k_1)$ is the unique solution of (\ref{3.3}) with $(\nu,d,\theta)=(\mu,1, 1)$, i.e., $k_1=k(\mu,1, 1)$.

Assume on the contrary that $g_\infty=\infty$. Let $(z_2,s_2)$ be the unique solution of (\ref{4.2}) with $\Gamma=(0,d,a,\beta,g_0)$ and $z_0(x)=u_0(x)$. Then $z_2(t,x)\ge u(t,x)$, $s_2(t)\ge g(t)$ for all $t\ge 0$ and $0\le x\le g(t)$. Similarly to the above, $s_2(t)-k(\beta,d,a) t\to \varsigma_2\in\mathbb{R}$ as $t\to\infty$.

Because of $(\beta,\mu)\in{\cal A}$, we have $k_1>k(\beta,d,a)$. This implies $s_1(t)-g(t)\ge s_1(t)-s_2(t)\to\infty$ and $\min_{0\le x\le g(t)}q_1(k_1 t+\varsigma_1-x)\to 1$ as $t\to\infty$. Thus, upon using (\ref{5.1}), $\lim_{t\to\infty}\min_{0\le x\le g(t)}z_1(t,x)=1$. There exists $t_2>t_1$ such that $z_1(t,x)>1-\ep$ for all $t\ge t_2$ and $0\le x\le g(t)$. Consequently, $v(t,x)>1-\ep$, and hence $a-u-bv<a-b(1-\ep)-u<0$ for all $t\ge t_2$ and $0\le x\le g(t)$. Take advantage of \cite[Lemma 3.2]{HW}, it follows that $g_\infty<\infty$. \ \ \ \fbox{}

\section{Estimates of asymptotic spreading speeds of $u,v$ and asymptotic speeds of $g,h$}
\setcounter{equation}{0} {\setlength\arraycolsep{2pt}

The authors of \cite{LinG} and \cite{Pan}, by means of the construction of the appropriate and elaborate upper and lower solutions, established some interesting results for the asymptotic spreading speeds of solution to the following Cauchy problem
 \bes\left\{\begin{array}{lll}
 u_t-du_{xx}=u(a-u-bv), &t>0,\ \ x\in\mathbb{R},\\[.1mm]
  v_t-v_{xx}=v(1-v+cu),\ \ &t>0, \ \ x\in\mathbb{R},\\[.1mm]
  u(0,x)=u_0(x), \ \ v(0,x)=v_0(x), &\ x\in\mathbb{R},\\[.1mm]
  0\le u_0(x)\le a, \ \ 0\le v_0(x)\le 1+ac, &\ x\in\mathbb{R}.
  \end{array}\right.\lbl{6.1}\ees
Their conclusions show that the prey and predator may have different asymptotic spreading speeds.

\begin{defi}\lbl{d6.1}Let $w(t,x)$ be a nonnegative function for $t>0$ and $x\in[0,\infty)$. A number $c_*>0$ is called the asymptotic spreading speed of $w(t,x)$ if

{\rm(a)} $\dd\lim_{t\to\infty}\sup_{x\ge(c_*+\ep)t}w(t,x)=0$ for any given $\ep>0$;

{\rm(b)} $\dd\lim_{t\to\infty}\inf_{0\le x\le(c_*-\ep)t}w(t,x)>0$ for any given $0<\ep<c_*$.
\end{defi}

The asymptotic spreading speed gives the observed phenomena imagining an observer moves to the right at a fixed speed \cite{WLL}, and it describes the speed at which the geographic range of the new population expands in population dynamics \cite{HsuZ}. Thus the asymptotic spreading of prey and predator are useful and important in understanding the interspecies action between the prey and predator. The background of prey-predator system implies that the predator has a negative effect on the prey, while the prey has a positive effect on the predator (see \cite{OL} for some biological results). Intuitively, we guess (believe) that the asymptotic propagation of prey (asymptotic spreading speed of $u$ and asymptotic speed of $g$) may be slower than the case of no predator, and that of the predator (asymptotic spreading speed of $v$ and asymptotic speed of $h$) may be faster than the case of no prey. However, our results indicate that this is not necessarily right.

The other related works on the asymptotic spreading speeds of evolutionary systems, please refer to  \cite{BHN, LWL, LZ, LL, Wh, WLL, Zhao} and the references cited therein. In some evolutionary systems, the nonexistence of constant asymptotic spreading speed has been observed, see Berestycki et al. \cite{Ber} for some examples.

In this section we study the asymptotic spreading speeds of
$u,v$ and asymptotic speeds of $g,h$. Assume $a>b(1+ac)$. In consideration of (\ref{3.5}), using the known results of (\ref{3.1}) and comparison principle, we see that both prey and predator must spread successfully as long as their moving parameters are suitably large. That is, there are $\beta_1,\,\mu_1>0$ such that $g_\infty=h_\infty=\infty$ for all $\beta\ge\beta_1$ and $\mu\ge\mu_1$.

Throughout this section we assume $a>b(1+ac)$, which is equivalent to $bc<1$ and $a>b/(1-bc)$. Denote
 \bess
 &c_1=2\sqrt{da}, \ \ c_2=2\sqrt{1+ac}, \ \ c_3=2\sqrt{da-db(1+ac)},& \\[.1mm]
 &c_4=2\sqrt{da-db},\ \ \ \ c_5=2\sqrt{(1+ac)(1-bc)}.&\eess

 \begin{theo}\lbl{th6.1} For any given $0<\ep\ll 1$, there exist $\beta_\ep,\,\mu_\ep,\,T\gg 1$ such that, when $\beta\ge\beta_\ep$ and $\mu\ge\mu_\ep$,
  \bes
  &u(t,x)=0\ \ \ {\rm for} \ t\ge T, \, \ x\ge(c_1+\ep)t,&\lbl{6.2a}\\[1mm]
  &v(t,x)=0\ \ \ {\rm for} \ t\ge T, \, \ x\ge(c_2+\ep)t, &\lbl{6.2}\\[1mm]
 &\dd\liminf_{t\to\infty}\min_{0\le x\le(c_3-\ep)t}u(t,x)\ge a-b(1+ac), \ \
\liminf_{t\to\infty}\min_{0\le x\le(2-\ep)t}v(t,x)\ge 1.&\lbl{6.3}\ees
\end{theo}

{\bf Proof}. According to the first limit of (\ref{3.4}), it follows that
 \bess
 \lim_{\beta\to\infty}\bar k_\beta=c_1, \ \ \lim_{\mu\to\infty}\bar k_\mu=c_2,\ \
 \lim_{\beta\to\infty}\ud k_\beta=c_3, \ \ \lim_{\mu\to\infty}\ud k_\mu=2, \eess
where $\bar k_\beta,\,\ud k_\mu,\,\bar k_\mu$ and $\ud k_\beta$ are given by (\ref{3.6})-(\ref{3.9}), respectively. Note that (\ref{3.6})-(\ref{3.9}), for any given $0<\ep\ll 1$, there exist $\beta_\ep\gg 1$ and $\mu_\ep\gg 1$ such that
 \bes
 &\dd c_3-\ep/2<\ud k_\beta\le\liminf_{t\to\infty}\frac{g(t)}t,\ \ \ \limsup_{t\to\infty}\frac{g(t)}t\le\bar k_\beta< c_1+\ep/4, \ \ \ \forall \ \beta\ge\beta_\ep,&\lbl{6.4} \\[1mm]
 &\dd 2-\ep/2<\ud k_\mu\le\liminf_{t\to\infty}\frac{h(t)}t, \ \ \ \limsup_{t\to\infty}\frac{h(t)}t\le\bar k_\mu< c_2+\ep/2, \ \ \ \forall \ \mu\ge\mu_\ep.&\lbl{6.5}\ees
As a conclusion, we can find a $\tau_1\gg 1$ such that, for all $t\ge \tau_1$, $\beta\ge\beta_\ep$ and $\mu\ge\mu_\ep$,
 \bes(c_3-\ep)t<g(t)<(c_1+\ep/2)t, \ \ \  (2-\ep)t<h(t)<(c_2+\ep)t.\lbl{6.6}\ees
Obviously, (\ref{6.2a}) and \eqref{6.2} hold.
Similarly to Step 1 in the proof of Theorem \ref{th4.4}, we can prove (\ref{6.3}). The proof is finished. \ \ \ \fbox{}

\begin{theo}\lbl{th6.2} Suppose $da<1$. Then the following hold:

{\rm(i)} For any given $0<k_0<c_3$, $(\ref{4d.2})$ holds as long as  $\beta$ and $\mu$ are suitably large.

{\rm(ii)} There exists $\mu_0\gg1$ such that, when $\mu>\mu_0$,
 \bes
  \limsup_{\beta\to\infty}\limsup_{t\to\infty}\frac{g(t)}t\le c_4.\lbl{6.12}\ees

{\rm(iii)} For any given $\ep>0$, there exist $\beta_\ep,\,\mu_\ep,\,T\gg 1$ such that, when $\beta\ge\beta_\ep$ and $\mu\ge\mu_\ep$,
 \bes
 &u(t,x)=0 \ \ \ {\rm for} \ \ t\ge T, \ \, x\ge(c_4+\ep)t,&\lbl{6.7} \\[.5mm]
 &\dd\lim_{t\to\infty}\sup_{x\ge(2+\ep)t}v(t,x)=0,&\lbl{6.8} \\[.5mm]
 &\dd\lim_{t\to\infty}\max_{(c_1+\ep)t\le x\le(2-\ep)t}|v(t,x)-1|=0.\lbl{6.9}& \ees

{\rm(iv)} There exists $\beta_0\gg1$ such that, when $\beta>\beta_0$,
 \bes
   \lim_{\mu\to\infty}\lim_{t\to\infty}\frac{h(t)}{t}=2.\lbl{6.10}\ees
\end{theo}

{\bf Proof}. Take advantage of \eqref{6.4} and \eqref{6.5}, we have
 \[\dd \liminf_{\beta\to\yy}\liminf_{t\to\infty}\frac{g(t)}t\ge c_3, \ \ \
 \liminf_{\mu\to\infty}\liminf_{t\to\infty}\frac{h(t)}t\ge 2.\]
Since $da<1$, it is obvious that
 \[2>c_1>c_4>c_3.\]
The conclusion (i) can be proved by the same way as that of Theorem \ref{th4.4}.

(ii) Choose $0<\ep\ll 1$ such that $c_1+\ep<2-\ep$. Then, in view of (\ref{6.6}), we have
 \bes
 g(t)<(c_1+\ep/2)t<(c_1+\ep)t<(2-\ep)t<h(t), \ \ \forall \ \beta\ge\beta_\ep, \ \mu\ge\mu_\ep, \ t\ge \tau_1.\lbl{6.11}\ees
The second conclusion of (\ref{6.3}) shows that for any given $0<\delta\ll 1$, there exists $\tau_2>\tau_1$ such that $v(t,x)\ge 1-\delta$ for all $t\ge \tau_2$ and $0\le x\le (2-\ep)t$. Combining this with (\ref{6.11}), we see that, when $\beta\ge\beta_\ep$ and $\mu\ge\mu_\ep$, $u$ satisfies
 \bess
 \left\{\begin{array}{lll}
 u_t-du_{xx}\le u[a-b(1-\delta)-u],\ \ &t\ge \tau_2,\ 0<x<g(t),\\[0.5mm]
 u_x(t,0)=0,\ \ u(t,g(t))=0,\ &t\ge \tau_2,\\[0.5mm]
 g'(t)=-\beta u_x(t,g(t)),\ &t\ge \tau_2.
 \end{array}\right.
 \eess
It follows that
 \[\limsup_{\beta\to\infty}\limsup_{t\to\infty}\frac{g(t)}t\le 2\sqrt{da-db(1-\delta)},\]
and then (\ref{6.12}) holds because $\delta>0$ is arbitrary.

(iii) The result (\ref{6.7}) is a direct consequence of (\ref{6.12}).

Now we prove (\ref{6.8}). By virtue of $c_4<2$ and (\ref{6.12}), there exist $\tau_*\gg 1$, $\beta_*\gg 1$ such that $g(t)<2t$ for all $t\ge \tau_*$ and $\beta\ge\beta_*$. This implies $u(t,x)=0$ for all $t\ge \tau_*,\,x\ge 2t$ and $\beta\ge\beta_*$.
Define
   \[s(t)=\max\{h(t),\,(2+\ep)t\} \ \ \ {\rm for}\ \ t\ge \tau_*.\]
Note that $v(t,x)=0$ for $x\ge h(t)$ and $v_x(t,h(t))<0$, it is not hard to see that $v$ satisfies, in the weak sense,
  \bess
 \left\{\begin{array}{lll}
  v_t-v_{xx}\le v(1-v),\ \ &t\ge \tau_*,\ \ 2t\le x<s(t),\\[0.5mm]
 v(t,x)\le M_2, \ &t\ge \tau_*,\ \ 2t\le x<s(t),\\[.5mm]
 v(t,s(t))=0, \ &t\ge \tau_*,
 \end{array}\right.
 \eess
where $M_2$ is given by Lemma \ref{l2.1}. Define
  \[\xi(t,x)=M_2{\rm e}^{s(\tau_*)-2\tau_*}{\rm e}^{2t-x},\ \ \ t\geq \tau_*, \ 2t\le x<s(t).\]
Clearly,
 \[\sup_{x\ge(2+\ep)t}\xi(t,x)\le M_2{\rm e}^{s(\tau_*)-2\tau_*}{\rm e}^{-\ep t}\to 0\]
as $t\to\infty$, and
 \bess
 &\xi(t,2t)>M_2, \ \xi(t,s(t))>0, \ \ t\ge \tau_*,&\\[.5mm]
  &\xi(\tau_*,x)> M_2, \ \ 2\tau_*\le x<s(\tau_*).&
  \eess
It is easy to verify that
  \[\xi_t-\xi_{xx}\ge \xi(1-\xi),\ \ \ t\geq \tau_*, \ 2t\le x<s(t).\]
By the comparison principle, $v(t,x)\le \xi(t,x)$ for all $t\ge \tau_*$ and $2t\le x<s(t)$. The limit (\ref{6.8}) is obtained.

In the following we prove (\ref{6.9}). Based on (\ref{6.11}), we see that $v$ satisfies
 \bess
 \left\{\begin{array}{lll}
  v_t-v_{xx}=v(1-v),\ \ &t\ge \tau_1,\, \ (c_1+\ep/2)t\le x<h(t),\\[0.5mm]
 v(t,(c_1+\ep/2)t)\le M_2, \ v(t,h(t))=0, \ \ &t\ge \tau_1,\\[.5mm]
 v(\tau_1,x)\le M_2, \ &(c_1+\ep/2)\tau_1\le x<h(\tau_1),
 \end{array}\right.
 \eess
where $M_2$ is given by Lemma \ref{l2.1}. Define
 \[\vp(t,x)=1+M_2{\rm e}^{h(\tau_1)}{\rm e}^{(c_1+\ep/2)t-x},\ \ \ t\ge \tau_1,\ (c_1+\ep/2)t\le x<h(t).\]
The direct calculations yield
 \[\vp_t-\vp_{xx}\ge \vp(1-\vp),\ \ \ t\ge \tau_1,\ (c_1+\ep/2)t\le x<h(t),\]
and
 \bess
  &\vp(t,(c_1+\ep/2)t)>1+M_2, \ \ \vp(t,h(t))\ge 1, \ \ \forall \ t\ge\tau_1,&\\[.2mm]
  &\vp(\tau_1,x)\ge 1+M_2, \ \ \forall \ (c_1+\ep/2)\tau_1\le x\le h(\tau_1).&
 \eess
By the comparison principle, $v(t,x)\le\vp(t,x)$ for $t\ge \tau_1$ and $(c_1+\ep/2)t\le x<h(t)$. According to (\ref{6.11}), we have $(c_1+\ep)t<h(t)$ for all $t\ge\tau_1$ and $\beta\ge\beta_\ep, \ \mu\ge\mu_\ep$. And so
 \bes\max_{x\ge (c_1+\ep)t}v(t,x)=\max_{(c_1+\ep)t\le x\le h(t)}v(t,x)\le\max_{(c_1+\ep)t\le x\le h(t)}\vp(t,x)=1+M_2{\rm e}^{h(\tau_1)}{\rm e}^{-\ep t/2}, \ \ \lbl{6.13}\ees
which implies $\dd\limsup_{t\to\infty}\max_{x\ge (c_1+\ep)t}v(t,x)\le 1$. This combined with the second inequality of (\ref{6.3}) allows us to derive (\ref{6.9}).

(iv) For any given $0<\sigma\ll 1$ and  $\beta\ge\beta_\ep$. Let $(q(y),k)$ be the unique solution of (\ref{3.3}) with $(\nu,d,\theta)=(\mu, 1, 1+\sigma)$. Then $q'(y)>0$,  $q(y)\to 1+\sigma$ as $y\to\infty$ and $\lim_{\mu\to\infty}k=2\sqrt{1+\sigma}$. Combining these facts with (\ref{6.11}) and (\ref{6.13}), we can find three constants $\mu_0>\mu_\ep,\,\tau_0>\tau_1,\,y_0\gg 1$ such that, for all $\mu\ge\mu_0$,
  \bess
  &k>c_1+\ep, \ \ \ h(t)>(c_1+\ep)t, \ \ \ \forall \ t\ge\tau_0,&\\[,3mm]
  &v(t,x)<(1+M_2{\rm e}^{h(\tau_1)}{\rm e}^{-\ep t/2})q(y),\ \ \ \forall \ t\ge\tau_0, \ x\ge (c_1+\ep)t, \ y\ge y_0.&\eess
Denote $K=M_2{\rm e}^{h(\tau_1)}$ and define
 \bess
 \bar h(t)&=&kt+\varrho K\big({\rm e}^{-\ep\tau_0/2}
 -{\rm e}^{-\ep t/2}\big)+y_0+h(\tau_0), \ \ t\ge\tau_0,\\[.3mm]
 \bar v(t,x)&=&\big(1+K{\rm e}^{-\ep t/2}\big)q(\bar h(t)-x), \ \ t\ge\tau_0, \ \
 (c_1+\ep)t\le x\le\bar h(t),
 \eess
where $\varrho$ is a positive constant to be determined. Obviously,
 \bess
 &\bar h(\tau_0)>h(\tau_0), \ \ \bar v(\tau_0,x)\ge v(\tau_0,x), \ \ \ \forall \ (c_1+\ep)\tau_0\le x\le h(\tau_0),&\\[.3mm]
 &\bar v(t,\bar h(t))=0, \ \ \bar v(t,(c_1+\ep)t)>v(t,(c_1+\ep)t), \ \ \
 \forall \ t\ge\tau_0.&\eess
In the same way as the arguments of \cite[Lemma 3.5]{KM} we can verify that, when $\varrho$ is suitably large,
 \bess
  &\bar v_t-\bar v_{xx}\ge\bar v(1-\bar v),\ \ \ \forall \ t\ge \tau_0,\, \ (c_1+\ep)t\le x<\bar h(t),&\\[0.5mm]
  &\bar h'(t)\ge-\mu\bar v_x(t,\bar h(t)), \ \ \ \forall \ t\ge \tau_0.&\eess
Because $v$ satisfies $v_t-v_{xx}=v(1-v)$ for $t\ge\tau_0$ and $(c_1+\ep)t\le x<h(t)$, by the comparison principle we have $v(t,x)\le\bar v(t,x)$ and $h(t)\le\bar h(t)$ for all $t\ge\tau_0$ and $(c_1+\ep)t\le x<h(t)$. Hence,
 \[\limsup_{t\to\infty}\frac{h(t)}t\leq k, \ \ \ \limsup_{\mu\to\infty}\limsup_{t\to\infty}\frac{h(t)}t\leq\limsup_{\mu\to\infty} k=2\sqrt{1+\sigma}.\]
The arbitrariness of $\sigma$ leads to
 \[\limsup_{\mu\to\infty}\limsup_{t\to\infty}\frac{h(t)}t\leq 2.\]
This together with the first inequality of (\ref{6.5}) derive (\ref{6.10}), and the proof is finished. \ \ \ \fbox{}

\begin{theo}\lbl{th6.3} If $d[a-b(1+ac)]>1+ac$, we have the following conclusions:

{\rm(i)} There exists $\beta_0\gg1$ such that, when $\beta>\beta_0$,
  \[\liminf_{\mu\to\infty}\liminf_{t\to\infty}\frac{h(t)}t\ge c_5.\]

{\rm(ii)} For any given $\ep>0$, there exist $\beta_\ep,\,\mu_\ep\gg 1$ such that
 \bes
\liminf_{t\to\infty}\min_{0\le x\le(c_1-\ep)t}u(t,x)>0, \ \ \ \
 \liminf_{t\to\infty}\min_{0\le x\le(c_5-\ep)t}v(t,x)\ge(1+ac)(1-bc)>0\lbl{6.14}\ees
 provided $\beta\ge\beta_\ep$, $\mu\ge\mu_\ep$.

{\rm(iii)} There exists $\mu_0\gg 1$ such that, when $\mu>\mu_0$,
 \bess
   \lim_{\beta\to\infty}\lim_{t\to\infty}\frac{g(t)}{t}=2\sqrt{da}=c_1.\eess

{\rm(iv)}  For any given $0<k_0<c_5$, $(\ref{4d.2})$ holds as long as $\beta$ and $\mu$ are suitably large.
 \end{theo}

{\bf Proof}. The assumption $d[a-b(1+ac)]>1+ac$ implies
$c_3>c_2$. Choose $\ep>0$ is so small that $c_3-\ep>c_2+\ep$, then, by (\ref{6.6}),
we have
 \bes
 h(t)<(c_2+\ep)t<(c_3-\ep)t<g(t), \ \ \forall \ \beta\ge\beta_\ep, \ \mu\ge\mu_\ep, \ t\ge \tau_1.\lbl{6.15}\ees

For any given $0<\delta\ll 1$,  there exists $0<\sigma_0\ll 1$ such that $2\sqrt{(1+ac)(1-bc)-c\sigma}>c_5-\delta$ for all $0<\sigma\le\sigma_0$, where $c_5=2\sqrt{(1+ac)(1-bc)}$. For such a fixed $\sigma$, combining (\ref{6.15}) with the first inequality of (\ref{6.3}), we have that there exists $\tau_3>\tau_1$ such that $(v,h)$ satisfies
  \bess
 \left\{\begin{array}{lll}
 v_t-v_{xx}\ge v\big[1+c(a-b(1+ac)-\sigma)-v\big], &t\ge \tau_3,\ \ 0<x<h(t),\\[0.5mm]
 v_x(t,0)=0,\ \ v(t,h(t))=0,\ &t\ge \tau_3,\\[0.5mm]
 h'(t)=-\mu v_x(t,h(t)),\ &t\ge \tau_3
 \end{array}\right.
 \eess
for all $\beta\ge\beta_\ep$, $\mu\ge\mu_\ep$. Let $(\ud v,\ud h)$ be the unique solution of \bess
 \left\{\begin{array}{lll}
 \ud{v}_t-\ud{v}_{xx}=\ud{v}\big[1+c(a-b(1+ac)-\sigma)-\ud{v}\big],\ \ &t\ge \tau_3,\ \ 0<x<\ud{h}(t),\\[0.5mm]
 \ud{v}_x(t,0)=0,\ \ \ud{v}(t,\ud{h}(t))=0,\ &t\ge \tau_3,\\[0.5mm]
 \ud{h}'(t)=-\mu \ud{v}_x(t,\ud{h}(t)),\ &t\ge \tau_3,\\[.5mm]
 \ud h(\tau_3)=h(\tau_3), \ \ \ud{v}(\tau_3,x)=v(\tau_3,x),\ &0\le x\le\ud h(\tau_3).
 \end{array}\right.
 \eess
Then $h(t)\ge\ud h(t),\,v(t,x)\ge\ud v(t,x)$ for all $t\ge\tau_3$ and $0\le x\le \ud h(t)$ by the comparison principle. Make use of   \cite[Theorem 3.1]{ZhaoW15}, it follows that, for any given $0<\ep\ll 1$,
 \bess
 &\dd\lim_{\mu\to\infty}\lim_{t\to\infty}\frac{\ud h(t)}t=2\sqrt{(1+ac)(1-bc)-c\sigma}>c_5-\delta,&\\[.5mm]
 &\dd\lim_{t\to\infty}\max_{0\le x\le(c_5-\delta-\ep)t}\big|\ud v(t,x)-[(1+ac)(1-bc)-c\sigma]\big|=0.&\eess
Consequently,
 \[\liminf_{\mu\to\infty}\liminf_{t\to\infty}\frac{h(t)}t\ge c_5, \ \ \ \liminf_{t\to\infty}\min_{0\le x\le(c_5-\ep)t}v(t,x)\ge(1+ac)(1-bc)\]
since $\delta$ and $\sigma$ are arbitrary. The conclusion (i) and the second inequality of (\ref{6.14}) are obtained.

Now we prove the first inequality of (\ref{6.14}). Notice $v\le M_2$ and the first inequality of (\ref{6.3}), there exists $\tau_4\gg 1$ such that $v(t,x)\leq\frac{2M_2}{a-b(1+ac)}u(t,x)$ for all $t\ge\tau_4$ and $0\le x\le(c_3-\ep)t$. In view of (\ref{6.15}) we see that $v(t,x)=0$ when $x\ge(c_3-\ep)t$. Denote $r=\frac{2M_2}{a-b(1+ac)}$, then $(u,g)$ satisfies
   \bess
 \left\{\begin{array}{lll}
 u_t-du_{xx}\ge u(a-u-bru), &t\ge \tau_4,\ \ 0<x<g(t),\\[0.5mm]
 u_x(t,0)=0,\ \ u(t,g(t))=0,\ &t\ge \tau_4,\\[0.5mm]
 g'(t)=-\beta u_x(t,g(t)),\ &t\ge \tau_4.
 \end{array}\right.
 \eess
Similarly to the above, we can show that
 \bes\liminf_{\beta\to\infty}\liminf_{t\to\infty}\frac{g(t)}t\ge 2\sqrt{ad}=c_1, \ \ \ \liminf_{t\to\infty}\min_{0\le x\le(c_1-\ep)t}u(t,x)\ge a/(1+br).\lbl{6.16}\ees
The above second inequality implies the first one of (\ref{6.14}).

The conclusion (iii) can be derived from the first inequality of (\ref{6.16}) and the second one of (\ref{6.4}).

Since $c_5<c_1$, the proof of (iv) is the same as that of Theorem \ref{th4.4}. The proof of Theorem \ref{th6.3} is complete. \ \ \ \fbox{}

\section{Discussion--biological significance of the conclusions}
\setcounter{equation}{0}

In this paper, we investigated a free boundary problem which describes the expanding of prey and predator in a one-dimensional habitat. In this model, the prey occupying the interval $[0,g(t)]$, while the predator with the territory $[0,h(t)]$ at time $t$. Here, the two free boundaries $x=g(t)$ and $x=h(t)$ may intersect each other as time evolves. They describe the spreading fronts of prey and predator, respectively. Our aim is to study its dynamics. Because these two free boundaries may intersect each other, it seems very difficult to understand the whole dynamics of this model.

{\bf (A)}\, Concerning the long time behaviors of solution, we established some realistic and more sophisticated results.

(I)\, If the prey (predator) species can not spread into $[0,\infty)$, then it will die out in the long run.

(II)\, When both two species spread successfully. For the weak predation case $b<a$, under the condition $a>b(1+ac)$, we find an important expanding phenomenon: If an observer were to move to
the right at a fixed speed less than $\min\{\ud k_\beta,\,\ud k_\mu\}$, it will be observed that
the two species will stabilize at the unique positive equilibrium state, while, if we observe the two species in front of the curves $x=g(t)$ and $x=h(t)$, we could see nothing. This is different from the Cauchy problem (\ref{6.1}) because in (\ref{6.1}) the two species become positive for all $x$ once $t$ is positive.

{\bf (B)}\, Main results about the spreading and vanishing show the following important phenomena, these look more realistic and may play an important role in the understanding of ecological complexity.

(I)\, When one of the initial habitat and moving parameter of the predator is ``suitably large", the predator is always able to successfully spread.

(II)\, If one of the initial habitat and moving parameter of the prey is ``suitably large", but both the initial habitat and moving parameter of the predator are ``suitably small", the former will spread successfully, while the latter will vanishes eventually.

(III)\,  When prey's initial habitat is much larger than that of predator, and predator spreads slowly, the prey will spread successfully and its territory always cover that of the predator, whether or not the latter spreads successfully.

(IV)\, In the case of strong predation, if the prey spreads slowly and the predator does quickly, the former will vanish eventually (it will be eaten up by the latter) and the predator will spread successfully.

{\bf (C)}\, The conclusions regarding the asymptotic propagations reveal the complicated and realistic spreading phenomena of prey and predator.

(I)\, When (\ref{1.1}) is uncoupled ($b=c=0$), the prey and predator satisfy
 \bes
 \left\{\begin{array}{lll}
 w_t-dw_{xx}=w(a-w), &t>0,\ \ 0<x<\gamma(t),\\[0.5mm]
  w_x(t,0)=w(t,\gamma(t))=0,\ \ \ &t\ge0,\\[0.5mm]
 g'_0(t)=-\beta w_x(t,\gamma(t)),\ \  &t\ge0,\\[0.5mm]
 w(0,x)=u_0(x),\ \ &0\le x\le \gamma=\gamma(0)\\[0.5mm]
  \end{array}\right.\lbl{1.1a}
 \ees
and
  \bes
 \left\{\begin{array}{lll}
 z_t-z_{xx}=z(1-z),\ \ &t>0, \ \ 0<x<\zeta(t),\\[0.5mm]
 v_x(t,0)=z(t,\zeta(t))=0,\ \ \ &t\ge0,\\[0.5mm]
 \zeta'(t)=-\mu v_x(t,\zeta(t)), \ &t\ge0,\\[0.5mm]
z(0,x)=v_0(x),\ \ & 0\le x\le h_0=\zeta(0),
 \end{array}\right.\lbl{1.1b}
 \ees
respectively. By use of \eqref{3.4},
 \[\lim_{\beta\to\yy}\lim_{t\to\yy}\frac{\gamma(t)}t=2\sqrt{da}, \ \ \
 \lim_{\mu\to\yy}\lim_{t\to\yy}\frac{\zeta(t)}t=2.\]
Similarly to the discussion in Section 6 we can show that
the asymptotic spreading speed of $w$ is $2\sqrt{da}$ when $\beta$ is sufficient large, and that of $z$ is $2$ when $\mu$ is sufficient large.

For the case that $a>b(1+ac)$ and $da<1$. The conclusions of Theorems \ref{th6.1} and \ref{th6.2} show that the asymptotic spreading speed of predator is $2$ and that of prey is between $2\sqrt{da-db(1+ac)}$ and $2\sqrt{da-db}$ when $\beta$ and $\mu$ are sufficiently large. Moreover,
 \bess
  &\dd\lim_{\mu\to\infty}\lim_{t\to\infty}\frac{h(t)}{t}=2 \ \ \ {\rm for} \ \
  \beta>\beta_0,&\\[1mm]
  &\dd 2\sqrt{da-db(1+ac)}\le\liminf_{\beta\to\infty}\liminf_{t\to\infty}\frac{g(t)}t
  \le \limsup_{\beta\to\infty}\limsup_{t\to\infty}\frac{g(t)}t\le 2\sqrt{da-db}.&\eess
These illustrate that the prey is not helpful to the predator's asymptotic propagation, while the predator could decrease that of the prey. The reason is that prey's ability to diffuse and grow is weaker than that of the predator.

In the case of $a>b(1+ac)$ and $d[a-b(1+ac)]\ge 1+ac$. The conclusions of Theorems \ref{th6.1} and \ref{th6.3} indicate that the asymptotic spreading speed of prey is  $2\sqrt{da}$ and that of predator is between $2\sqrt{(1+ac)(1-bc)}$ and $2\sqrt{(1+ac)}$  when $\beta$ and $\mu$ are sufficiently large. Moreover,
 \bess
  &\dd\lim_{\beta\to\infty}\lim_{t\to\infty}\frac{g(t)}{t}=2\sqrt{da} \ \ \ {\rm for} \ \
  \mu>\mu_0,&\\[1mm]
  &\dd 2\sqrt{(1+ac)(1-bc)}\le\liminf_{\mu\to\infty}\liminf_{t\to\infty}\frac{h(t)}t
  \le \limsup_{\mu\to\infty}\limsup_{t\to\infty}\frac{h(t)}t\le 2\sqrt{1+ac}.&\eess
Since $2\sqrt{(1+ac)(1-bc)}>2$ in this case, we see that the prey accelerates the asymptotic propagation of predator, while the predator has no effect on that of the prey. The reason is that the prey spreads faster and provides predator with more food.

(II)\, In the case of $a>b(1+ac)$ and $da<1$. If an observer were to move to the right at a fixed speed less than $2\sqrt{da-db(1+ac)}$, it will be observed that
the two species will stabilize at the unique positive equilibrium state; When the observer do this with a fixed speed $k\in(2\sqrt{da-db},\, 2)$, he can only watch the predator; When the observer do this with a fixed speed  greater than $2$, he can not find anything because the two species have not arrived in his horizon.

(III)\, For the case $a>b(1+ac)$ and $d[a-b(1+ac)]> 1+ac$. When we  are to move to the right at a fixed speed less than $2\sqrt{(1+ac)(1-bc)}$, we shall observe that the two species will stabilize at the unique positive equilibrium state; When we do this with a fixed speed $k\in(2\sqrt{(1+ac)},\,2\sqrt{da})$, we can only see the prey; When we do this with a fixed speed greater than $2\sqrt{da}$, we could see nothing because the two species are not in our sight.

{\bf (D)}\, A great deal of previous mathematical investigation on the spreading of population has been based on the traveling wave fronts of prey-predator system (\ref{6.1}). A striking difference between our free boundary problem (\ref{1.1}) and the Cauchy problem (\ref{6.1}) is that the spreading fronts in (\ref{1.1}) are given explicitly by two functions $x=g(t)$ and $x=h(t)$, beyond them respectively the population densities of prey and predator are zero, while in (\ref{6.1}), the two species become positive for all $x$ once $t$ is positive. Secondly, (\ref{6.1}) guarantees successful spreading of the two species for any nontrivial initial populations $u(0, x)$ and $v(0,x)$, regardless of their initial sizes and supporting area, but the dynamics of (\ref{1.1}) possesses the multiplicity and complexity  of spreading and vanishing. The phenomena exhibited by these multiplicities and complexities seem closer to the reality.

\section*{\bf Appendix. Proof of \eqref{x1}}
 \def\theequation{A.\arabic{equation}}\setcounter{equation}{0}

Set $W=w_1-w_2$, $G=g_1-g_2$, then $W, G$ satisfy \eqref{2.14}, \eqref{2.15} and the estimate \eqref{b.10} holds. Using (\ref{2.15}) and \eqref{b.10} we have
  \bes\|G'\|_{C^{\frac\alpha 2}([0,T])}&\leq& \beta\|g^{-1}_1W_y\|_{C^{\frac\alpha 2,0}(\Delta_T^1)}
 +\beta\|(g_1^{-1}-g_2^{-1})w_{2y}\|_{C^{\frac\alpha 2,0}(\Delta_T^1)}\nonumber\\[.5mm]
 &\leq& C_{11}\big(\|W\|_{C^{\frac{1+\alpha} 2,\,1+\alpha}(\Delta_T^1)}+\|G\|_{C^1([0,T])}\big)\nonumber\\[.5mm]
 &\leq& C_{12}\big(\|z_1-z_2\|_{C(\Delta_T^1)}+\|G\|_{C^1([0,T])}\big).\lbl{A.1}\ees

Recall $W(0,y)=0$, $G(0)=G'(0)=0$. Take advantage of the mean value theorem and \eqref{b.10}, \eqref{A.1}, it follows that
 \bess
 \|W\|_{C(\Delta_T^1)}
 &\leq&T^{\frac\alpha 2}\|W\|_{C^{\frac\alpha 2,0}(\Delta_T^1)}\leq C_4T^{\frac\alpha 2}\big(\|z_1-z_2\|_{C(\Delta_T^1)}+\|G\|_{C^1([0,T])}\big),\\[1mm]
 \|G\|_{C^1([0,T])}&\leq& 2T^{\frac\alpha 2}\|G'\|_{C^{\frac\alpha 2}([0,T])}
 \leq 2C_{12}T^{\frac\alpha 2}\big(\|z_1-z_2\|_{C(\Delta_T^1)}+\|G\|_{C^1([0,T])}\big).\eess
Thus we have \bess
  \|W\|_{C(\Delta_T^1)}+\|G\|_{C^1([0,T])}\leq (C_4+2C_{12})T^{\frac\alpha 2}\big(\|z_1-z_2\|_{C(\Delta_T^1)}+\|G\|_{C^1([0,T])}\big).
  \eess
If we choose $T>0$ such that $(C_4+2C_{12})T^{\frac\alpha 2}\le 1/2$, then
  \bess
  \|W\|_{C(\Delta_T^1)}+\|G\|_{C^1([0,T])}\leq \|z_1-z_2\|_{C(\Delta_T^1)},
  \eess
which is exactly \eqref{x1}.

 \vskip 4pt \noindent {\bf Acknowledgment:} The author would like to thank the anonymous referees for their helpful comments and suggestions.

\end{document}